\documentclass[12pt]{article}
\usepackage{amsmath,amsthm,amscd,amssymb}
\usepackage{latexsym}
\begin{document}

\numberwithin{equation}{section}

\title{Isoperimetric and Universal Inequalities\\for Eigenvalues}

\author{Mark S.\ Ashbaugh\thanks{Partially supported by National
Science Foundation (USA) grants DMS-9500968 and DMS-9870156.}\\
Department of Mathematics\\University of Missouri\\Columbia, MO
65211-0001\\e-mail: mark@math.missouri.edu}

\maketitle

\begin{description}\item[1991 Mathematics Subject Classification:]
Primary 35P15, Secondary 58G25, 49Rxx.

\item[Keywords and phrases:] eigenvalues of the Laplacian,
isoperimetric inequalities for eigenvalues, Faber-Krahn inequality,
Szeg\H{o}-Weinberger inequality, Payne-P\'{o}lya-Weinberger conjecture,
Sperner's inequality, biharmonic operator, bi-Laplacian, clamped
plate problem, Rayleigh's conjecture, buckling problem, the
P\'{o}lya-Szeg\H{o} conjecture, universal inequalities for eigenvalues,
Hile-Protter inequality, H.~C.\ Yang's inequality.

\item[Short title:] Isoperimetric and Universal
Inequalities\end{description}

\begin{abstract}  This paper reviews many of the known inequalities
for the eigenvalues of the Laplacian and bi-Laplacian on bounded
domains in Euclidean space.  In particular, we focus on isoperimetric
inequalities for the low eigenvalues of the Dirichlet and Neumann
Laplacians and of the vibrating clamped plate problem (i.e., the
biharmonic operator with ``Dirichlet'' boundary conditions).  We also
discuss the known universal inequalities for the eigenvalues of the
Dirichlet Laplacian and the vibrating clamped plate and buckling
problems and go on to present some new ones.  Some of the names
associated with these inequalities are Rayleigh, Faber-Krahn,
Szeg\H{o}-Weinberger, Payne-P\'{o}lya-Weinberger, Sperner, Hile-Protter,
and H.~C.\ Yang.  Occasionally, we will also comment on extensions of
some of our inequalities to bounded domains in other spaces, 
specifically, $S^n$ or $H^n$.
\end{abstract}

\section{Introduction}

\subsection{The Eigenvalue Problems}

The first eigenvalue problem we shall introduce is that of the {\it
fixed membrane}, or {\it Dirichlet Laplacian}.  We consider the
eigenvalues and eigenfunctions of $-\Delta$ on a bounded domain
(=connected open set) $\Omega$ in Euclidean space ${\Bbb R}^n$,
i.e., the problem
\begin{align}\label{eq1.1a} \tag{1.1a} -\Delta u&=\lambda u\quad
\text{in }
\Omega \subset {\Bbb R}^n,\\
\label{eq1.1b}\tag{1.1b} u&=0\quad \text{on } \partial \Omega.\end{align} 
It is well-known that this problem has a real and
purely discrete spectrum $\{\lambda_i\}_{i=1}^\infty$ where
\setcounter{equation}{1}
\begin{equation}\label{eq1.2} 0< \lambda _1<\lambda_2\leq
\lambda_3\leq
\ldots \to \infty .\end{equation} 
Here each eigenvalue is repeated
according to its multiplicity.  An associated orthonormal basis of
real eigenfunctions will be denoted $u_1$, $u_2$, $u_3$, $\ldots$.  
In fact, throughout this paper we will assume that all functions we 
consider are real.  This entails no loss of generality in the present 
context.

The next problem we introduce is that of the {\it free membrane}, or
{\it Neumann Laplacian}.  This is the problem
\begin{align} \tag{1.3a}\label{eq1.3} -\Delta v&=\mu v\quad
\text{in }\Omega \subset {\Bbb R}^n,\\
\tag{1.3b} \frac{\partial v}{\partial n}&=0\quad
\text{on }\partial \Omega.\end{align} 
Here $\partial /\partial n$ denotes the outward normal derivative 
on $\partial \Omega$, where we now assume that $\partial \Omega$ 
is sufficiently smooth.  With this assumption, problem (1.3) has 
spectrum $\{ \mu_i\}_{i=0}^\infty$ where
\setcounter{equation}{3}\begin{equation} \label{eq1.4} 0=\mu_0<
\mu_1\leq \mu_2\leq
\ldots \to \infty ,\end{equation} with the eigenvalues again
repeated according to their multiplicities.  A corresponding
orthonormal basis of real eigenfunctions will be denoted $\{
v_i\}_{i=0}^\infty$.

Next we introduce the {\it clamped plate problem}, or eigenvalue
problem for the {\it Dirichlet biharmonic operator} (for an explanation 
of this terminology see \cite{Dav2} or \cite{Owen}), which describes
the characteristic vibrations of a clamped plate.  This problem is
given by
\begin{align} 
\label{eq1.5a} \tag{1.5a} 
\Delta ^2w&=\Gamma w\quad
\text{in }\Omega \subset {\Bbb R}^n,\\
\label{eq1.5b} \tag{1.5b} w&=0=\frac{\partial w}{\partial n}\quad
\text{ on } \partial \Omega .
\end{align} 
We will denote the eigenvalues and an associated orthonormal basis 
of real eigenfunctions by $\{ \Gamma _i\}_{i=1}^\infty$ and
$\{w_i\}_{i=1}^\infty$, respectively.  The eigenvalues
$\Gamma_i$ satisfy 
\setcounter{equation}{5}\begin{equation}\label{eq1.6}0<\Gamma_1\leq
\Gamma_2\leq
\Gamma_3\leq
\ldots \to \infty.\end{equation}

Lastly, we introduce the {\it buckling problem}, which determines
the critical buckling load of a clamped plate subjected to a uniform
compressive force around its boundary (this description applies to 
the $n=2$ case of the problem).  This problem again 
involves the biharmonic operator and is formulated as
\begin{align}\tag{1.7a}\label{eq1.7a} \Delta^2v&=-\Lambda \Delta
v\quad \text{in } \Omega \subset {\Bbb R}^n,\\
\tag{1.7b}\label{eq1.7b}v&=0=\frac{\partial v}{\partial n}\quad
\text{on } \partial \Omega .\end{align} 
It also has a discrete spectrum consisting of positive eigenvalues 
of finite multiplicity with infinity as their only accumulation point.  
We denote the eigenvalues by $\{ \Lambda _i\}_{i=1}^\infty$ and a 
corresponding orthonormal basis of real eigenfunctions by 
$\{ v_i\}_{i=1}^\infty$.  Thus
\setcounter{equation}{7} \begin{equation}\label{eq1.8} 0<\Lambda_1\leq
\Lambda_2\leq \Lambda _3\leq \ldots \to \infty .\end{equation}

Good sources of information on many of these problems include the 
books of Bandle \cite{Ban}, B\'{e}rard \cite{Ber}, Courant-Hilbert 
\cite{CH}, Davies \cite{Dav}, Kesavan \cite{Kes2}, Leis \cite{Leis}, 
P\'{o}lya and Szeg\H{o} \cite{PS}, Reed and Simon \cite{ReSi}, and 
Safarov and Vassiliev \cite{SV}.  The review papers of Davies 
\cite{Dav2}, Payne \cite{Pa2}, \cite{Pa3}, \cite{Pa4}, Protter 
\cite{Prot}, and Talenti \cite{Ta4}, \cite{Ta5} are also quite 
useful, as is Leissa's monograph \cite{Leissa} on the vibration 
of plates.  Rayleigh's classic {\it The Theory of Sound} 
\cite{Ray} is highly recommended as collateral reading.

\subsection{Rearrangement}

In this subsection we introduce the notion of {\it spherically 
symmetric rearrangement} (or {\it Schwarz symmetrization}) and 
some of its properties.  Suppose that we have a bounded measurable 
function $f$ on the bounded measurable set $\Omega \subset 
{\Bbb R}^n$.  We can consider its {\it distribution function} 
$\mu_f(t)$ defined by
\begin{equation}\label{eq1.9}\mu_f(t)=|\{ x\in \Omega ||f(x)|>
t\}|\end{equation} where $|\cdot |$ denotes Lebesgue measure.  The
distribution function can be viewed as a function from
$[0,\infty)$ to $[0,|\Omega |]$.  It is clearly a nonincreasing
function.  The {\it decreasing rearrangement of
$f$}, denoted $f^*$, is essentially the inverse of $\mu_f$ and is 
defined by
\begin{equation}\label{eq1.10}f^*(s)=\inf \{ t\geq
0|\mu_f(t)<s\}.\end{equation} 
It is a nonincreasing function on $[0,|\Omega |]$.

Before defining the spherically decreasing rearrangement of a
function we define the {\it spherical} (or {\it symmetric}) 
{\it rearrangement} of a set.  For a
bounded measurable set $\Omega \subset {\Bbb R}^n$, we define its
spherical rearrangement $\Omega ^\star$ as the ball centered at the
origin having the same measure as $\Omega $, i.e.,
$|\Omega ^\star|=|\Omega |$.  We can now define the {\it spherically 
(symmetric) decreasing rearrangement} $f^\star :\Omega^\star \to 
{\Bbb R}$ by
\begin{equation}\label{eq1.11} f^\star (x)=f^*(C_n|x|^n)\quad
\text{for}\quad x\in \Omega^\star \end{equation} where
$C_n=\pi^{n/2}/\Gamma \left(\frac{n}{2}+1\right) =$ volume of the
unit ball in ${\Bbb R}^n$.  A verbal description of $f^\star $ runs as
follows: $f^\star $ is that function of $x\in \Omega ^\star $ which is
spherically symmetric, radially decreasing (in the weak sense of
nonincreasing) and equimeasurable with $|f|$, i.e., $f^\star$ and $|f|$
share the same distribution function (see (1.9) above).  One way to view 
this is that for
each of the level sets $\Omega_t\equiv\{x\in \Omega ||f(x)|>t\}$ of $|f|$
we take its spherical rearrangement $\Omega^\star_t$ and define this to
be the corresponding level set of $f^\star$ (which, we recall, is placed
concentrically with $\Omega^\star$).  This viewpoint gives an
interpretation of $f^\star$ which avoids the intermediate function $f^*$
of a single (=volume) variable.  Good sources of further information
on rearrangements are \cite{Ban}, \cite{GMGT}, \cite{Gun}, 
\cite{HLP}, \cite{Kaw}, \cite{Lieb1}, \cite{Lieb2}, \cite{LiLo}, 
\cite{PS}, \cite{Simon}, \cite{Ta0}, \cite{Ta1}, \cite{Ta2}, 
\cite{Ta4}, \cite{Ta5}.

A key fact, which is evident from the equimeasurability
of $f$, $f^*$, and $f^\star $, is that
\begin{equation}\label{eq1.12} \int_\Omega
f^2=\int^{|\Omega|}_0f^*(s)^2ds=\int_{\Omega^\star}(f^\star)^2.
\end{equation}
There is no reason for not putting $|f|$, $f^*$, and $f^\star$ as
integrands here, or even $|f|^p$, $|f^*|^p$, and
$|f^\star |^p$, except that \eqref{eq1.12} is all we need later in this
paper.

We shall be particularly concerned with how spherical rearrangement
affects the first eigenfunction $u_1$ of the Dirichlet Laplacian.
The key property here is that spherical rearrangement typically
decreases, and in any case cannot increase, the Dirichlet norm $\int
_\Omega |\nabla u_1|^2$.  In particular, it is known that for any
function $f$ in the Sobolev space $H^1_0(\Omega )$ (for background
and notation on Sobolev spaces we recommend
\cite{Ad}, \cite{Brez}, \cite{Kes2}, and \cite{LiLo}) $f^\star 
\in H^1_0(\Omega ^\star)$ and
\begin{equation}\label{eq1.13} \int_{\Omega^\star} |\nabla f^\star |^2\leq
\int_\Omega |\nabla f|^2.\end{equation} This inequality is crucial
to the proof of the Faber-Krahn inequality given in Section 2.  Note,
too, that by
\eqref{eq1.12} the $L^2$-norm of $u_1$ does not change when we
replace it by $u^\star_1$.  For a discussion of
\eqref{eq1.13}, see Glaser, Martin, Grosse, and Thirring \cite{GMGT}, 
Gunson \cite{Gun}, Kawohl \cite{Kaw}, Lieb \cite{Lieb1}, 
\cite{Lieb2}, Lieb and Loss \cite{LiLo}, P\'{o}lya-Szeg\H{o} \cite{PS}, 
and/or Talenti \cite{Ta0}, \cite{Ta1},
\cite{Ta4}, \cite{Ta5}.

We also note a certain elementary property of rearrangement as it
affects integrals of products of functions.  This is that for two
nonnegative measurable functions $f$ and $g$ on $\Omega$
\begin{equation}\label{eq1.14} \int_\Omega fg\leq
\int_{\Omega^\star }f^\star g^\star .\end{equation} 
There is a corresponding lower bound if we rearrange $f$ and $g$ in 
opposite senses.  For this we need the notion of {\it spherically 
(symmetric) increasing rearrangement}, which we denote by a lower 
$\star $.  The definition is almost identical to that of
spherically decreasing rearrangement, except that $g_\star $ should be
radially increasing (in the weak sense) on $\Omega^\star $.  Then we
have, for example, 
\begin{equation}\label{eq1.15} \int_\Omega fg\geq \int_{\Omega^\star }
f^\star g_\star .\end{equation}

In fact, these inequalities take their most elementary form if we use 
``signed'' rearrangements of $f$ and $g$, that is, we define $f^\star$, 
$g^\star$, $f_\star$, $g_\star$, etc., in terms of the distribution 
function as given in \eqref{eq1.9} {\it except} that we use just $f(x)$ 
in place of $|f(x)|$ in that formula.  With $f^\star$, $g^\star$, etc., 
defined in this alternative way, \eqref{eq1.14} and \eqref{eq1.15} 
hold without any need to restrict $f$ and $g$ to be nonnegative 
functions.  All of this basic material is admirably presented in 
\cite{HLP}, where the essential features of the process (similarly 
ordered, oppositely ordered) are brought to the fore.  In a sense 
the essence of the whole business is the simple algebraic inequality
\begin{equation}\label{eq1.16} (a-b)(c-d)\geq 0\end{equation}
if $(a,b)$ and $(c,d)$ are similarly ordered (and the reverse if these 
vectors are oppositely ordered).  Thus, if $\vec{v}=(a,b)$ and 
$\vec{w}=(c,d)$ are similarly ordered
\begin{equation}\label{eq1.17} \vec{v}\cdot \vec{w}=ac+bd\geq ad+bc=\vec{v}
\cdot (d,c).\end{equation}
The relevance of this simple vector inequality for \eqref{eq1.14} (in its 
general ``signed'' version) is that we can approach \eqref{eq1.14} by 
approximating $f$ and $g$ by simple functions decomposed over sets of equal 
measures.  In this setting the relevance of \eqref{eq1.17} is apparent 
(as a first approximation in the single-variable setting one might think 
of passing to Riemann sums over subintervals of equal length).  See 
\cite{HLP}, and the more recent article by Baernstein \cite{Baern}, 
for further development of these ideas.  

We come finally to a result on how rearrangement affects the
solution to the Poisson equation on a bounded domain $\Omega$ with 
homogeneous Dirichlet boundary conditions
\begin{align}\label{eq1.18a}\tag{1.18a}-\Delta u&=f \quad \text{in }
\Omega \subset {\Bbb R}^n,\\
\label{eq1.18b}\tag{1.18b}u&=0\quad \text{on } \partial \Omega .\end{align} 
Suppose we solve this problem for $u$ and then take its
spherically decreasing rearrangement $u^\star $.  We want to compare 
$u^\star$ to the solution $v$ of the symmetrized (=spherically 
rearranged) problem 
\begin{align} \tag{1.19a}\label{eq1.19a} -\Delta v&=f^\star\quad
\text{in }\Omega^\star ,\\
\tag{1.19b}\label{eq1.19b} v&=0 \quad \text{on } \partial
\Omega^\star.\end{align} It turns out that if $u\geq 0$ on 
$\Omega $, then
\setcounter{equation}{19}\begin{equation} \label{eq1.20} 
0\leq u^\star \leq
v.\end{equation} 
For $f^\star $ in \eqref{eq1.19a} we can even use a signed 
rearrangement of $f$ (still under the assumption that $f$ is 
sufficiently positive that $u\geq 0$ on $\Omega$, which implies, 
in particular, that $\int_\Omega f\geq 0$).  For example, if $f$ 
is bounded below by $c$ we can take $f^\star =(f-c)^\star +c$; 
this gives a radially decreasing function that passes into 
negative values out near $\partial \Omega $ if $f$ itself is 
not always nonnegative.

The key to proving \eqref{eq1.20} is a differential inequality for
$u^\star $ (or equivalently $u^*$) and the fact that (1.19)
reduces to a one-dimensional problem that can be solved explicitly.
In particular, we have (in the radial variable $r=|x|$) 
\begin{equation}\label{eq1.21} -\frac{1}{r^{n-1}}
(r^{n-1}v')'=f^\star \end{equation} and hence
\begin{equation} \label{eq1.22}
r^{n-1}v'(r)=-\int^r_0\tau^{n-1}f^\star (\tau )d\tau \end{equation} 
and finally
\begin{equation}\label{eq1.23}
v(r)=\int^R_rt^{1-n}\int^t_0\tau^{n-1}f^\star(\tau ) d\tau
dt,\end{equation} 
where we have incorporated the conditions $v'(0)=0$ (necessary for 
$v$ to be smooth at the origin in $\Omega^\star$) 
and $v(R)=0$, with $R$ denoting the radius of $\Omega ^\star$.
Recalling the condition $\int_\Omega f\geq 0$, we see that this is
precisely the condition needed to keep $v$ nonnegative, since this
guarantees that $\int^t_0\tau^{n-1}f^\star (\tau )d\tau \geq 0$ for 
all $t\in [0,R]$ (recall, too, that up to the constant factor $n C_n$,
$\int^t_0\tau^{n-1}f^\star(\tau )d\tau $ is 
$\int ^s_0f^*(\sigma )d\sigma$ where $\sigma$ is now the volume 
variable $\sigma =C_n\tau^n$ and similarly $s$ is related to $t$ 
by $s=C_n t^n$; in general, we shall reserve $s$ and $\sigma$ as 
``volume variables'', while $t$, $\tau$, $r$, and $\rho$ will 
be used as radial variables).

The basic idea used to get \eqref{eq1.20} goes back at least to
Weinberger \cite{Wein2}.  The method is also given by Bandle in her
book \cite{Ban} and figures prominently in works of Talenti
\cite{Ta0}, \cite{Ta1}, \cite{Ta2}, and Chiti \cite{Chiti1}, 
\cite{Chiti2}, \cite{Chiti3}, \cite{Chiti4}.  The extended form 
presented above (allowing a signed rearrangement of $f$)
occurs in Talenti \cite{Ta3}, where it is instrumental in his
treatment of the first eigenvalue $\Gamma_1$ of the clamped plate
problem.  As such, it also figures in the later papers on the subject
by Nadirashvili \cite{Nad1},
\cite{Nad2}, \cite{Nad3} and Ashbaugh and Benguria \cite{AB9} (see also
\cite{ABL}).  A useful discussion of the extended form is also found
in the papers of Kesavan \cite{Kes1}, \cite{Kes3}, \cite{Kes4}, 
\cite{Kes5}.  The interaction between rearrangements and partial 
differential equations encompasses a variety of topics by a wide 
range of authors.  A small sampling of other works in this field 
includes \cite{ADLT}, \cite{ALT1}, \cite{ALT2}, \cite{ALT3}, 
\cite{ATDL}, and \cite{Baern}.   

That one can allow the function $f$ to have variable sign while using a 
signed rearrangement of it so long as the condition $u \geq 0$ holds 
was perhaps first observed by Talenti.  Earlier authors had assumed 
$f \geq 0$, apparently because this is a clean and easily checked 
condition from which $u \geq 0$ follows immediately via the maximum 
principle.  At any rate, throughout this paper when we invoke the result
\eqref{eq1.20} under the condition $u \geq 0$ while employing a signed
rearrangement of $f$, we shall refer to it as Talenti's theorem.   

To show \eqref{eq1.20}, one applies a fairly standard procedure to
\eqref{eq1.18a} integrated over level sets of
$u$.  One then uses the classical isoperimetric inequality and
arrives at an integral-differential inequality for
$u^*$ in which $v'$ (for $v$ as given by \eqref{eq1.23}) appears 
on one side.  Inequality \eqref{eq1.20} then follows by integration.  
For details, one might consult the papers of Talenti or Kesavan 
mentioned above, or books such as \cite{Ban} or \cite{Kaw}.  

\subsection{The Rayleigh-Ritz Inequality}

Throughout this paper we shall have many occasions to use the
Rayleigh-Ritz inequality, which gives a simple way to bound
eigenvalues from above based on trial functions.  For example, for
the Dirichlet Laplacian on the bounded domain $\Omega $ one has
\begin{equation}\label{eq1.24} \lambda_1(\Omega ) =\inf_{{\varphi
\in D(-\Delta )\atop \varphi \not\equiv 0}}
\frac{\int_{\Omega }\varphi (-\Delta \varphi )}{\int_\Omega
\varphi^2}\end{equation} 
where $\varphi$ is a real trial function in
the domain of $-\Delta $ (denoted $D(-\Delta )$).  One can also get
at the higher eigenvalues by imposing orthogonality conditions on
the class of trial functions used.  For example,
\begin{equation}\label{eq1.25} \lambda_{k+1}(\Omega)=\inf_{{\varphi
\in D(-\Delta )\atop 0\not\equiv \varphi \bot u_1,\ldots
,u_k}}\frac{\int_\Omega \varphi (-\Delta \varphi)}{\int_\Omega
\varphi^2}\end{equation} 
where $u_1,\ldots, u_k$ denote the first $k$
eigenfunctions of $-\Delta$ on $\Omega$.  Beyond this, and somewhat
more useful for our purposes below, there is a quadratic form
formulation of the Rayleigh-Ritz inequality.  For the Dirichlet
Laplacian it reads as follows:
\begin{equation} \label{eq1.26}\lambda_{k+1}(\Omega )=\inf_{{\varphi
\in H^1_0(\Omega )\atop 0\not\equiv \varphi \bot u_1,\ldots ,u_k}}
\frac{\int_\Omega |\nabla \varphi |^2}{\int_\Omega
\varphi^2}.\end{equation} 
We note that the Sobolev space
$H^1_0(\Omega )=Q(-\Delta )=$ the form domain of $-\Delta$ in this
case (see, for example \cite{Dav}, \cite{Kes2}, and \cite{ReSi}).  This
formulation has the advantage over \eqref{eq1.25} that the trial
function $\varphi$ can be chosen from a larger class of functions
(and, in particular, it need have essentially only one
square-integrable derivative, not two).

The Rayleigh-Ritz inequality (and the closely related Min-Max Principle) 
applies to any semi-bounded (from below) self-adjoint operator on a
Hilbert space.  For more details and discussion, the reader might 
consult Bandle \cite{Ban}, B\'{e}rard \cite{Ber}, Chavel \cite{Cha1},
\cite{Cha2}, Davies \cite{Dav}, Kesavan
\cite{Kes2}, and/or Reed and Simon, vol.~4 \cite{ReSi}.

\section{Isoperimetric Inequalities for Eigenvalues}

\subsection{The Faber-Krahn Inequality}

One of the earliest isoperimetric inequalities for an eigenvalue is
certainly that for the first eigenvalue of the Dirichlet Laplacian
(= fixed membrane problem) conjectured by Rayleigh \cite{Ray} in
1877:
\begin{equation}\label{eq2.1} \lambda_1(\Omega )\geq
\lambda_1(\Omega ^\star )\quad \text{for } \Omega \subset {\Bbb
R}^n\end{equation} with equality if and only if $\Omega $ is a ball,
i.e., $\Omega =\Omega^\star $.  This result was subsequently proved
(independently) by Faber \cite{Fab} and Krahn \cite{Kr1},
\cite{Kr2} in the 1920's using symmetrization.  In terms of spherical
symmetrization (=spherical rearrangement) the proof can now be
reduced to a few lines.  One uses the first eigenfunction $u_1$ for
$\Omega $ and \eqref{eq1.12}, \eqref{eq1.13} to conclude
\begin{equation}\begin{split}\label{eq2.2} \lambda_1(\Omega ) & =
\int_\Omega |\nabla u_1|^2\\ &\geq \int_{\Omega ^\star }|\nabla u_1^*|^2\\
&\geq \lambda_1(\Omega ^\star ),\end{split}\end{equation} 
where the last line follows from the Rayleigh-Ritz inequality for 
$\lambda_1$ of $-\Delta $ on $\Omega ^\star $ and the fact (mentioned 
in connection with \eqref{eq1.13}) that $u_1\in H^1_0(\Omega )$ 
implies that $u^\star _1\in H ^1_0(\Omega^\star )$, and thus that 
$u^\star _1$ is a valid trial function in the Rayleigh-Ritz inequality 
for $\lambda_1(\Omega^\star )$.  The characterization of the case of 
equality is somewhat technical, so we refer the reader to the 
literature for it.  Good sources are Kawohl's book \cite{Kaw}, 
or Kesavan's articles \cite{Kes3}, \cite{Kes4}.

\subsection{The Szeg\H{o}-Weinberger Inequality}

We next turn to the isoperimetric result for $\mu_1(\Omega )$, the
first nonzero Neumann eigenvalue, due to Szeg\H{o} \cite{Sz2} (for
$n=2$ and $\Omega$ simply connected) and Weinberger \cite{Wein1} (in
full generality).  This is the next simplest (or even perhaps the 
simplest) isoperimetric result for eigenvalues.  The result was 
first suggested by Kornhauser and Stakgold \cite{KoSt}, who also 
obtained some results in support of it.  The Szeg\H{o}-Weinberger 
result states that
\begin{equation}\label{eq2.3} \mu_1(\Omega )
\leq \mu_1(\Omega^\star)\quad\text{for } 
\Omega \subset {\Bbb R}^n,\end{equation} with
equality if and only if $\Omega $ is a ball.

We follow Weinberger's method of proof, since it is very natural and
lends itself to the problem we treat next, that of finding the
optimal upper bound for $\lambda_2/\lambda_1$.  First we recall that
$\mu_1(\Omega )$ may be characterized through the variational principle
\begin{equation}\label{eq2.4} \mu_1(\Omega ) =\min_{{\varphi \in
H^1(\Omega )\atop 0\not\equiv \varphi \bot 1}}
\frac{\int_\Omega |\nabla \varphi |^2}{\int_\Omega
\varphi^2}.\end{equation} 
This is just the (quadratic form version of the) Rayleigh-Ritz 
inequality for $\mu_1$, since $v_0$ is a constant.  Following 
Weinberger, we take as trial functions
$\varphi=P_i$, $i=1,\ldots , n$, such that $\int_\Omega P_i=0$ for
$i=1,\ldots , n$ (this is proved via a topological argument, but
given \eqref{eq2.5} below all it says is that we can choose to place
our origin of coordinates at an appropriate generalized center of
mass) with
\begin{equation}\label{eq2.5} P_i(x)=g(r)\frac{x_i}{r},\end{equation}
where the $x_i$'s are Cartesian coordinates, $x=(x_1,\ldots ,x_n)\in
{\Bbb R}^n$, $r=|x|$, and
\begin{equation}\label{eq2.6} g(r)=\begin{cases} w(r)=\text{ 
[``right'' radial function for a ball } B_R \text{ of radius } R \\
\quad \quad \quad \quad 
\text{ where } B_R=\Omega^\star] 
\quad \quad \text{ for } 
0\leq r\leq R,\\ w(R) 
\qquad \qquad \qquad \quad \quad \quad \quad \quad 
\text{ for } r\geq R.\end{cases}
\end{equation} 
The function $w$ arises as
a solution of the radial equation when one separates variables on
the ball $B_R$ and is therefore basically a Bessel function.  We use
only the facts that $w(0)=0$ and that $w$ satisfies
\begin{equation}\label{eq2.7} w''+\frac{n-1}{r} w'-\frac{n-1}{r^2}
w+\mu _1(B_R)w=0\quad\text{for}\quad 0<r<R.\end{equation} 
Since $w'(R)=0$ and $\mu_1(B_R)$ is defined as the first eigenvalue 
of this boundary value problem, it follows (by an appropriate choice 
of sign) that $w(r)$ is increasing on $[0,R]$ and hence that $g$ 
is everywhere nondecreasing for $r\geq 0$.  By substituting our 
trial functions $P_i$ into the Rayleigh-Ritz inequality for 
$\mu_1$, we find
\begin{equation}\label{eq2.8} \mu_1(\Omega )\int_\Omega P^2_i\leq
\int_\Omega |\nabla P_i|^2.\end{equation} 
Summing this in $i$ for $1\leq i\leq n$, we arrive at
\begin{equation}\begin{split}\label{eq2.9} \mu_1(\Omega )&\leq
\frac{\int_\Omega \sum^n_{i=1} |\nabla P_i|^2}{\int_\Omega
\sum^n_{i=1}P^2_i} = \frac{\int_\Omega \left[g'(r)^2+\frac{n-1}{r^2}
g(r)^2\right]}{\int_\Omega g(r)^2}\\ &= \frac{\int_\Omega
B(r)}{\int_\Omega g(r)^2}\end{split}\end{equation} 
where we have defined
\begin{equation}\label{eq2.10} B(r)\equiv g'(r)^2+\frac{n-1}{r^2}
g(r)^2.\end{equation} 
Now $B(r)$ is easily seen to be decreasing for
$0\leq r\leq R$ by differentiating and using the differential
equation \eqref{eq2.7}.  One finds
\begin{equation}\label{eq2.11}
B'(r)=-2[\mu_1(B_R)gg'+(n-1)(rg'-g)^2/r^3]<0\quad\text{for}\quad
0<r<R.\end{equation} 
In addition, $B(r)=(n-1)w(R)^2/r^2$ for $r\geq R$
shows that $B$ is decreasing for $r\geq R$.  Since $B$ is 
continuous for all $r\geq 0$, it is clearly also decreasing there.
Now one has only to observe that
\begin{equation}\label{eq2.12} \int_\Omega B(r)\leq
\int_{\Omega^\star}B(r)\end{equation} since the volumes integrated over
are the same in both cases, while in passing from the left- to
right-hand sides you are exchanging integrating over $\Omega
\backslash \Omega^\star$ for integrating over $\Omega^\star\backslash
\Omega$ (which are sets of equal volume).  Since $B$ is (strictly) 
decreasing this clearly increases the value of
the integral unless $\Omega=\Omega^\star$, when equality obtains.
Similarly we find that
\begin{equation}\label{eq2.13} \int_\Omega g(r)^2\geq
\int_{\Omega^\star}g(r)^2\end{equation} since $g$ is nondecreasing. 
Thus we arrive at
\begin{equation}\label{eq2.14} \mu_1(\Omega )
\leq \frac{\int_{\Omega^\star}
B(r)}{\int_{\Omega^\star}g(r)^2}=\mu_1(\Omega^\star),\end{equation} 
since each $P_i$ is precisely an eigenfunction of $-\Delta$ with 
eigenvalue $\mu_1(B_R)$ for the domain $B_R=\Omega^\star$.  This 
completes the proof of the Szeg\H{o}-Weinberger inequality, including 
the characterization of the case of equality.

\subsection{The Payne-P\'{o}lya-Weinberger Conjecture}

The next isoperimetric result that we consider is that for
$\lambda_2/\lambda_1$ for the fixed membrane problem.  In 1955 and
1956, Payne, P\'{o}lya, and Weinberger \cite{PPW1}, \cite{PPW2} proved
that
\begin{equation}\label{eq2.15} \frac{\lambda_2}{\lambda_1}\leq
3\quad \text{for}\quad \Omega \subset {\Bbb R}^2\end{equation} and
conjectured that
\begin{equation}\label{eq2.16} \left.\frac{\lambda_2}{\lambda_1}\leq
\frac{\lambda_2}{\lambda_1}\right|_{\text{disk}}=
\frac{j^2_{1,1}}{j^2_{0,1}}\approx 2.5387 \end{equation}
with equality if and only if $\Omega$ is a disk (i.e., $\Omega
=\Omega^\star$) and where $j_{p,k}$ denotes the
$k^{\text{th}}$ positive zero of the Bessel function $J_p(t)$ 
(we follow the notation of Abramowitz and Stegun \cite{AS} here). 
For general dimension $n\geq 2$ the analogous statements are
\begin{equation}\label{eq2.17} \frac{\lambda_2}{\lambda_1}\leq
1+\frac{4}{n}\quad \text{ for } \Omega \in {\Bbb R}^n,\end{equation} and
the {\it PPW conjecture}
\begin{equation}\label{eq2.18} \left.\frac{\lambda_2}{\lambda_1}\leq
\frac{\lambda_2}{\lambda_1}\right|_{n-\text{ball}}=
\frac{j^2_{n/2,1}}{j^2_{n/2-1,1}},\end{equation}
with equality if and only if $\Omega$ is an $n$-ball.  This PPW
conjecture was proved in 1990 by Rafael Benguria and the author (see
\cite{AB1}, \cite{AB2}, \cite{AB4}).

We proceed now with the proof of \eqref{eq2.18}.  This proof follows
the main outline of Weinberger's proof of the $\mu_1$ result given
previously, but it is substantially more complicated due mainly to
the fact that, unlike $v_0=1/\sqrt{|\Omega |}$, $u_1$ is not
constant.  We start from the variational principle for $\lambda_2$
\begin{equation}\label{eq2.19} \lambda_2(\Omega )=\min_{{\varphi\in
H^1_0(\Omega)\atop 0\not\equiv \varphi \bot u_1}}
\frac{\int_{\Omega }|\nabla \varphi |^2}{\int_\Omega
\varphi^2},\end{equation} or, better for our purposes here,
\begin{equation}\label{eq2.20} \lambda_2(\Omega )-\lambda_1(\Omega
)\leq
\min_{\begin{array}{c}{\scriptstyle P u_1 \in H^1_0(\Omega )}\\
{\scriptstyle\int_\Omega P u^2_1 =0}\\ {\scriptstyle P\not\equiv
0}\end{array}} \frac{\int_\Omega |\nabla P|^2u^2_1}{\int_\Omega
P^2u^2_1},\end{equation} which follows from \eqref{eq2.19} via
integration by parts.  (Note that with $\varphi =Pu_1$, 
\begin{equation*}\begin{split} \int_\Omega |\nabla \varphi |^2
&=\int_\Omega [|\nabla
P|^2u^2_1+2Pu_1\nabla P\cdot \nabla u_1+P^2|\nabla u_1|^2]\\
&= \int_\Omega |\nabla P|^2u^2_1+
\int_\Omega [\nabla (P^2u_1)]\cdot \nabla u_1\\
&=\int_\Omega |\nabla P|^2u_1^2+\int_\Omega P^2u_1(-\Delta u_1)\\
&=\int_\Omega |\nabla P|^2u^2_1+
\lambda_1(\Omega )\int_\Omega P^2u^2_1,\end{split}\end{equation*}
where we integrated by parts in the second-to-last step.)  In
\eqref{eq2.20} we shall use $n$ trial functions
$P=P_i$, $i=1,\ldots ,n$, such that $\int_\Omega P_iu^2_1=0$ for
$i=1,\ldots ,n$ (again proved by a topological argument and
interpretable as a generalized center of mass result) where
\begin{equation}\label{eq2.21} P_i=g(r)\frac{x_i}{r}\end{equation}
and
\begin{equation}\label{eq2.22} g(r)=\begin{cases} w(r)=\text{
[``right'' radial function for a ball }B_R\\
\quad \quad \quad \quad \text{ of radius } R] 
\qquad \quad \text{ for } 0\leq r\leq R,\\ w(R) 
\qquad \qquad \quad \quad \quad \quad \quad \ \quad 
\text{ for } r\geq R.\end{cases}\end{equation} 
The right $R$ in this case turns out to be the unique $R$ such that 
$\lambda_1(B_R)=\lambda _1(\Omega)$.  This is a key fact and explains 
why at bottom our proof of the PPW conjecture is a ``fixed-$\lambda_1$ 
result''.  A major motivation for this is the comparison result 
of Chiti, to be explained shortly.  The right function $w(r)$ for 
$0\leq r\leq R$ is fairly complicated, being a ratio of Bessel 
functions.

From here we can proceed much as before finding
\begin{equation}\label{eq2.23} \lambda_2 (\Omega
)-\lambda_1(\Omega )\leq \frac{\int_\Omega B(r) u^2_1}{\int_\Omega
g(r)^2u^2_1}\end{equation} where
\begin{equation}\label{eq2.24} B(r)\equiv g'(r)^2+\frac{n-1}{r^2}
g(r)^2.\end{equation} 
Again it can be confirmed (though it is
harder this time around) that $B$ is decreasing (and positive) and
$g$ is increasing (and positive).  For the proof of these facts, see
\cite{AB1}, \cite{AB2}, \cite{AB4}, or
\cite{AB7}.  The proof found in \cite{AB1} (for $n=2$ and $3$) and
\cite{AB2} (for all $n\geq 2$) is based on the product
representation for the Bessel functions involved and ultimately
comes down to certain inequalities between Bessel function zeros
(see also Section 4 of \cite{AB4} for a brief further discussion of
these issues).  On the other hand, the proof found in \cite{AB4} (see
also \cite{AB7}) is somewhat simpler and lends itself to
generalization to a version of the PPW conjecture for domains in
$S^n$.  These topics are dealt with further in \cite{AB11}, \cite{AB12},
and also to some extent below.

To continue from \eqref{eq2.23} we use rearrangement as follows:
\begin{equation}\label{eq2.25} \lambda_2 (\Omega )-\lambda_1(\Omega
)\leq
\frac{\int_{\Omega^\star}B(r)^\star u^{\star 2}_1}{\int_{\Omega^\star }
g(r)^2_\star u_1^{\star 2}}\leq
\frac{\int_{\Omega^\star}B(r)u_1^{\star 2}}{\int_{\Omega^\star}g(r)^2
u_1^{\star 2}}.\end{equation}
The basic ideas here are that the integral of a product is increased if 
we similarly rearrange both
functions (see \eqref{eq1.14}), while it is decreased if we oppositely 
rearrange them (see \eqref{eq1.15}).
Note that in \eqref{eq2.25} we similarly rearrange in the
numerator, while we oppositely rearrange in the denominator.
The removal of the $\star$'s from $g(r)^2$ and $B(r)$ in the last step 
is allowed due to their respective monotonicity properties.

Finally, we need to invoke a comparison result due to Chiti
\cite{Chiti1}, \cite{Chiti2}, \cite{Chiti3},
\cite{Chiti4} to replace the $u_1^\star$'s in \eqref{eq2.25} by 
something more tractable.  This result says that if we take $z(r)$ as 
the normalized first eigenfunction of $B_R$ (recall that $R$ is so
chosen that $\lambda_1(\Omega )=\lambda_1(B_R)$) then $B_R \subset 
\Omega^\star$ and $z$ will be larger than $u_1^\star$ (thought of as 
a  function of $r$) for
$r$ near $0$, say for $0<r<r_1$, and smaller farther out, in this
case for $r>r_1$.  Under this condition one finds that if $u^\star_1$ 
is replaced by $z$ in the final member of \eqref{eq2.25} the numerator
can only go up and the denominator can only go down.  Thus
\begin{equation}\label{eq2.26} \lambda_2(\Omega )-\lambda_1(\Omega
)\leq
\frac{\int_{B_R}B(r)z^2}{\int_{B_R}g(r)^2z^2}=\lambda_2(B_R)-
\lambda_1(B_R),\end{equation}
the final equality holding because our choices of $g$ and $w$ were
made precisely to make things come back together in this way.  Since
$\lambda_1(B_R)=\lambda_1(\Omega )$, \eqref{eq2.26} implies
\begin{equation}\label{eq2.27} \lambda_2(\Omega )\leq
\lambda_2(B_R)\end{equation} and hence also
\begin{equation}\label{eq2.28} \frac{\lambda_2(\Omega
)}{\lambda_1(\Omega )}\leq
\frac{\lambda_2(B_R)}{\lambda_1(B_R)}=\left.\frac{\lambda_2}{\lambda_1}
\right|_{\text{any $n$-ball}}
=\frac{\lambda_2(\Omega^\star)}{\lambda_1(\Omega^\star)}.\end{equation} 
We note that the PPW inequality \eqref{eq2.28} was really proved as a
subsidiary result to \eqref{eq2.27}, which should be regarded as a
``fixed-$\lambda_1$ result''.  It is this result that is in a certain
sense more fundamental and which is most easily extended to other
settings, for example, domains contained in a hemisphere of $S^n$.

Indeed, in subsequent work \cite{AB7}, \cite{AB11}, \cite{AB12} we
have established the fixed-$\lambda_1$ result
\eqref{eq2.27} for domains contained in a hemisphere of $S^n$.  In
addition, from this result we can then derive the result
\begin{equation}\label{eq2.29} \frac{\lambda_2(\Omega
)}{\lambda_1(\Omega )}\leq
\frac{\lambda_2(\Omega^\star)}{\lambda_1(\Omega^\star)}\end{equation} 
for a domain $\Omega$ contained in a hemisphere of $S^n$.  Here
$\Omega^\star$ denotes the geodesic ball (=polar  cap) having the same
volume as $\Omega$ and $B_R$ must be taken as the (unique) polar cap
having the same
$\lambda_1$ as $\Omega$.  Furthermore, there is a
corresponding result for $\mu_1$ (that is, $\mu_1(\Omega )\leq
\mu_1(\Omega^\star)$) for a domain $\Omega$ contained in a hemisphere
of $S^n$ and for $\Omega$ a bounded domain in $H^n$.  See
\cite{AB7}, \cite{AB8} for more on this aspect of the $\mu_1$
problem (in particular, the result for $H^n$ is due to Chavel, while 
both Chavel and Bandle had earlier variants of the $S^n$ result). 
The Faber-Krahn inequality also extends in sharp form to
the spaces $S^n$ and $H^n$.  See, for example, Sperner \cite{Sp} or
Friedland and Hayman \cite{FrHa}.  In particular, Sperner's bound
$\lambda_1(\Omega) \geq \lambda_1(\Omega^\star)$ for $\Omega 
\subset S^n$ is one of the elements necessary to our proof of 
\eqref{eq2.27} and \eqref{eq2.29} for domains contained in a 
hemisphere of $S^n$.    

\subsection{Rayleigh's Conjecture for the Vibrating Clamped Plate}

We turn now to Rayleigh's conjecture for the vibration of a clamped
plate.  Rayleigh made this conjecture in 1877 in the first edition of
his book {\bf The Theory of Sound} \cite{Ray} (see p.~382 of volume
1 of the second edition).  In terms of the notation introduced in
Section 1, Rayleigh's conjecture states that
\begin{equation}\label{eq2.30} \Gamma_1 (\Omega )\geq \Gamma_1 
(\Omega^\star)\quad \text{for } \Omega
\subset {\Bbb R}^2,\end{equation} 
with equality if and only if $\Omega $ is itself a disk.  It is
natural to conjecture that this inequality might apply equally 
well in ${\Bbb R}^n$.

This conjecture seems to have lain dormant until around 1950, when
Szeg\H{o} made some progress on it (see
\cite{Sz1}, \cite{Sz3}, and also the treatment in P\'{o}lya and
Szeg\H{o}'s book \cite{PS}).  In fact, Szeg\H{o} was able to prove
\eqref{eq2.30} for simply connected domains $\Omega$ having a
nonnegative first eigenfunction,
$w_1\geq 0$.  At the time it seems to have been thought possible that
$w_1\geq 0$ for all domains.  However, it soon developed that this
cannot be expected to hold in general.  Results of Duffin \cite{Du}
and Duffin and Shaffer \cite{DuSh} were enough to disabuse people
of the notion that $w_1\geq 0$ should always hold (see also the much 
more recent article of Kozlov, Kondrat'ev, and Maz'ya \cite{KKM}, 
as well as further references given in \cite{ABL} and \cite{AL}).

The next advance toward the proof of Rayleigh's conjecture
\eqref{eq2.30} came in 1981 when Talenti \cite{Ta3} developed an
approach to the sign-indeterminate case using separate
rearrangements on the sets $\Omega_+=\{x\in
\Omega |w_1>0\}$ and $\Omega_-=\{x\in \Omega |w <0\}$.
This procedure led Talenti to two radial subproblems (tied to the
two separate balls $(\Omega_+)^\star$ and $(\Omega_-)^\star$ and 
also to the full ball $\Omega^\star$) which he
decoupled and considered as variational problems in their own right.
By this means, he obtained a lower bound to $\Gamma_1(\Omega )$
depending on the parameter $t\equiv |\Omega_+|/|\Omega|$ for $t\in
[0,1]$.  If this problem had been minimized at $t=0$ (and therefore also at
$t=1$, by symmetry) then Talenti's approach would have proved Rayleigh's
conjecture.  But unfortunately Talenti's minimum occurred at $t=1/2$
(for each $n=2,3,\ldots$) and he was only able to obtain inequalities
\begin{equation}\label{eq2.31} \Gamma _1(\Omega )\geq
d'_n\Gamma_1(\Omega^\star)\quad \text{for } \Omega \subset{\Bbb
R}^n.\end{equation} These bounds fall short of Rayleigh's
conjecture, since the constants $d'_n$ are less than $1$.  However,
Talenti was able to show that $d'_n\in (1/2,1)$, and, in particular,
he found $d'_2\approx 0.9777$, $d'_2\approx 0.7391$,
$d'_4\approx 0.6524$, \ldots .  In fact, the $d'_n$'s seem to
decrease monotonically to $1/2$ as $n$ goes to infinity (that
$\lim_{n\to \infty} d'_n=1/2$ is proved in \cite{AL}).

In 1992 Nadirashvili (see \cite{Nad1}, \cite{Nad2}, \cite{Nad3}) saw
how to improve upon the approach of Talenti by using slightly
different radial subproblems, each living only on
$(\Omega _+)^\star$ or $(\Omega_-)^\star$, but coupled together 
via boundary conditions.  With his approach and a geometric 
rearrangement argument, Nadirashvili was able to give a proof of
\eqref{eq2.30} for $n=2$.  However, it was not clear that the same
approach could handle cases with $n\geq 3$.

To that end, Rafael Benguria and the author introduced a more
analytical variant of Nadirashvili's approach, wherein the final
analysis comes down to understanding a certain function defined
explicitly in terms of Bessel functions.  (Bessel functions, while
certainly in the background in the clamped plate problem, were not
in evidence in any part of Nadirashvili's proof.)  By studying the
behavior of this function we were able to show \cite{AB9} that for
$n=2$ and 3 the minimizer of our parametrized two-ball variational
problem occurs when $t=|\Omega_+|/|\Omega|$ is 0 (or 1),
thereby proving Rayleigh's conjecture in dimensions 2 and 3.  In
subsequent work with Richard Laugesen \cite{AL} (see also
\cite{ABL}) it was shown that for all $n\geq 4$ things go the other
way, i.e., the parametrized two-ball minimizer occurs at $t=1/2$.
This yields bounds of the form
\begin{equation}\label{eq2.32} \Gamma_1 (\Omega ) > d_n\Gamma_1
(\Omega^\star)\quad \text{for } \Omega \subset {\Bbb R}^n, n\geq
4\end{equation} where our dimension-dependent constants $d_n$ turn
out always to be better than the corresponding $d'_n$'s of Talenti,
and, in fact, go to 1 as $n$ goes to infinity.  For example, we have
$d_4\approx 0.9537$,
$d_5\approx 0.9218$, $d_6\approx 0.9077$, \ldots ($d_8\approx
0.8998$ appears to be the minimum over all dimensions
$n$).

We now give the proof of Rayleigh's conjecture \eqref{eq2.30} for
dimensions $2$ and $3$.  As a warm-up, we first treat the case where
$w_1\geq 0$ on $\Omega$.  This is the case first successfully handled
by Szeg\H{o} \cite{Sz1}, \cite{Sz3} (see also P\'{o}lya-Szeg\H{o} 
\cite{PS} and Talenti \cite{Ta3}).  We recall the variational 
principle for $\Gamma_1$
\begin{equation}\label{eq2.33} \Gamma_1(\Omega )=\min_{{\varphi \in
H^2_0(\Omega )\atop \varphi \not\equiv 0}}
\frac{\int_\Omega (\Delta \varphi)^2}{\int_\Omega \varphi
^2}.\end{equation} In particular, this gives equality if we take
$\varphi=w_1$.  If we now decreasing rearrange $f=-\Delta w_1$ to
$f^\star=(-\Delta w_1)^\star$ where our rearrangement respects signs, 
we can invoke Talenti's rearrangement result \eqref{eq1.20} for the 
problem $-\Delta w_1=f$ in $\Omega $, $w_1=0$ on $\partial \Omega$, 
to get
\begin{equation}\label{eq2.34} \Gamma_1(\Omega )=\frac{\int_\Omega
(\Delta w_1)^2}{\int_\Omega w^2_1}\geq
\frac{\int_{\Omega^\star}[(-\Delta
w_1)^\star ]^2}{\int_{\Omega^\star}v^2}=\frac{\int_{\Omega^\star}(\Delta
v)^2}{\int_{\Omega^\star}v^2}\geq \Gamma_1(\Omega^\star),\end{equation}
which is the result we want, assuming that $v$ is truly an
admissible trial function for the Rayleigh-Ritz inequality (=
variational principle) for $\Gamma_{1}$ on $\Omega^\star$.  To verify
this, we note that $v$ is a radial $C^2$ function on the ball
$\Omega^\star$ and satisfies $v(R)=0$ where $R$ is the radius of 
$\Omega^\star$.  To see that $v$ also satisfies
$\frac{\partial v}{\partial n}=0$ on $\partial \Omega^\star$ we compute
(since $\frac{\partial w_1}{\partial n}=0$ on $\partial \Omega $)
\begin{equation}\label{eq2.35} 0=\int_{\partial \Omega
}\frac{\partial w_1}{\partial n}=-\int_\Omega (-\Delta
w_1)=-\int_{\Omega^\star} f^\star=-\int_{\Omega^\star}(-\Delta v)=
\int_{\partial
\Omega^\star}\frac{\partial v}{\partial n}\end{equation} 
and this shows that $\frac{\partial v}{\partial n}=0$ on 
$\partial \Omega^\star$ since $v$ is radial and
therefore $\int_{\partial
\Omega^\star}\frac{\partial v}{\partial n}=|\partial
\Omega^\star|v'(R)=n C_n R^{n-1} v'(R)$.  This proves Szeg\H{o}'s clamped
plate result, i.e., Rayleigh's conjecture for the clamped plate in
the case of a first eigenfunction of fixed sign.  This proof holds
for all dimensions $n$.

We now proceed to the proof of the general case for $n=2$ and $3$.  It
turns out that the proof does not work for $n\geq 4$ (though
Rayleigh's conjecture may well still be true there).  As mentioned
earlier, in the general case we proceed by decomposing the problem
into two (coupled) subproblems, each on a separate ball (on
$(\Omega_+)^\star $ and $(\Omega_-)^\star $ in the notation 
introduced above).  Following this decomposition procedure, we have

\begin{equation}\begin{split}\label{eq2.36} \Gamma_1(\Omega )
=\frac{\int_\Omega (\Delta w_1)^2}{\int_\Omega
w^2_1}&=\frac{\int_{\Omega^+}(\Delta w_1)^2+\int_{\Omega^-}(\Delta
w_1)^2}{\int_{\Omega^+}w^2_1+\int_{\Omega^-}w^2_1}\\ &\geq
\frac{\int_{B_a}(\Delta u)^2+\int_{B_b}(\Delta
v)^2}{\int_{B_a}u^2+\int_{B_b}v^2}\end{split}\end{equation} 
where
$B_a=$ ball of radius $a=(\Omega _+)^\star$ and $B_b=$ball of radius
$b=(\Omega_-)^\star$ and $u$ and $v$ satisfy symmetrized Dirichlet
problems on $B_a$ and $B_b$, respectively.  Note that the numerator
in \eqref{eq2.36} is unchanged by the introduction of $u$ and $v$,
while the denominator is (typically) increased (and certainly does
not decrease) by virtue of inequality
\eqref{eq1.20} as applied to $u$ and $v$ (and $(w^+_1)^\star $ and
$(w^-_1)^\star$), respectively.  Now, not only do $u$ and $v$ satisfy
\begin{equation}\label{eq2.37} u(a)=0=v(b)\end{equation} 
but $u'(a)$ and $v'(b)$ are coupled by the equation
\begin{equation}\label{eq2.38}
a^{n-1}u'(a)=b^{n-1}v'(b).\end{equation} 
This comes about by an argument similar to that in \eqref{eq2.35}. 
We have
\begin{equation}\begin{split}\label{eq2.39}0&=\int_{\partial \Omega}
\frac{\partial w_1}{\partial n}=\int_\Omega \Delta
w_1=-\int_{B_a}(-\Delta w_1)^\star +\int_{B_b}(\Delta w_1)^\star \\
&=\int_{B_a}\Delta u-\int_{B_b}\Delta v=\int_{\partial
B_a}\frac{\partial u}{\partial n}-\int_{\partial B_b}\frac{\partial
v}{\partial n}\\ &=
n C_n[a^{n-1}u'(a)-b^{n-1}v'(b)].\end{split}\end{equation} 
The seemingly odd changes of sign in certain terms come about because
we introduce a sign change for $\Omega_-$ so that we will still be
comparing positive functions (and thus $u$ and $v$ are both positive).  
As a curiosity we mention that in fact
\eqref{eq2.39} (and also \eqref{eq2.35}) uses only the fact that 
$\partial w_1 / \partial n$ has average value $0$ on the boundary of 
$\Omega$.  Thus our results all hold under the weaker assumption that 
our plate has its edge fixed and ``clamped on average'', a somewhat 
curious result first observed by Richard Laugesen.  

To finish the argument we view the final member of \eqref{eq2.36} as
a variational problem in its own right where we now treat $u$ and
$v$ as trial functions subject to \eqref{eq2.37},
\eqref{eq2.38}, and the condition that they be radial functions on
$B_a$ and $B_b$, respectively.  Thus we consider the minimization
problem
\begin{equation}\label{eq2.40} J(a)\equiv \inf_{\varphi ,\psi}
\frac{\int_{B_a}(\Delta \varphi )^2+\int_{B_b}(\Delta\psi
)^2}{\int_{B_a}\varphi^2+\int_{B_b}\psi ^2}\end{equation} with
$\varphi$ and $\psi$ radial and satisfying the boundary conditions
\eqref{eq2.37} and
\eqref{eq2.38} (with $\varphi$ and $\psi$ replacing $u$ and $v$ 
respectively).  This is what we refer to as our parametrized two-ball 
variational problem.  Note that $a$ now appears as a final parameter 
(which we will eventually need to minimize over, since we have no 
control on it).  We observe that, due to scaling freedom, we can assume 
that $a^n+b^n=1$, which is equivalent to normalizing $|\Omega|$ at $C_n$.
In \eqref{eq2.40}, and henceforth, we regard $b$ as defined
implicitly in terms of $a$ by this relation.  We also remark that in
terms of our earlier parameter $t=|\Omega_+|/|\Omega |$ we have
$t=a^n$.  In a sense $t$ is the preferred variable here, since as a
function of $t\in [0,1]\; J$ is symmetric about $1/2$.  We compromise 
by speaking in terms of $a^n$ at times as we proceed.

Leaving $a\in (0,1)$ fixed for the time being, standard
variational theory can be applied to \eqref{eq2.40} with the result
that we find the Euler equations
\begin{align}\tag{2.41a}\label{eq2.41a} \Delta^2\varphi &=\mu
\varphi \quad \text{in } B_a,\\
\tag{2.41b}\label{eq2.41b} \Delta^2\psi &=\mu \psi \quad \text{in }
B_b,\end{align} and boundary conditions
\begin{align}\tag{2.41c}\label{eq2.41c} \varphi (a)=0=\psi (b),\\
\tag{2.41d}\label{eq2.41d}
a^{n-1}\varphi'(a)=b^{n-1}\psi'(b),\end{align} and
\begin{equation}\tag{2.41e}\label{eq2.41e} \Delta \varphi (a)+\Delta
\psi (b)=0.\end{equation} The expressions in the last formula make
sense because $\varphi$ and $\psi$ are radial.  This formula comes
about as a product of the variational theory, and is what is known
as a ``natural'' boundary condition (see, for example, \cite{Gou} or
\cite{WeSt}).  See \cite{AB9} for more on its derivation and for more
details on this material in general (see also \cite{ABL} for a
general overview of the topic).

The boundary value problem (2.41) can now be solved
explicitly in terms of Bessel functions and modified Bessel
functions.  The result is an implicit relation for the eigenvalues
(and in particular the first eigenvalue $\mu_1=J(a)$) which we can
analyze to determine how
$J(a)$ behaves for $0\leq a\leq 1$.  Since 
\setcounter{equation}{41}\begin{equation}\label{eq2.42} \Gamma_1(\Omega ) 
\geq \min_aJ(a)\end{equation} 
and $J(0)=J(0^+)=J(1)=J(1^-)=\Gamma (\Omega^\star)$ (this
requires some work to see, though it certainly should already appear
plausible; essentially one wants to confirm that taking $a=0$ gives
the same result as the limit $a\to 0^+$ and similarly for $a=1$), we
will be done if we can show that $\min J(a)$ occurs at $a=0$ (and
$a=1$).  For dimensions 2 and 3 this can be shown (see \cite{AB9})
and this completes the proof of the general Rayleigh conjecture for
the clamped plate in these cases.

For $n\geq 4$, the analysis given above holds all the way until the
final step, where it is found that $\min_a J(a)$ occurs not at $a=0$
and 1, but at $a^n=1/2$.  Following this outcome to its logical
conclusion leads to the nonoptimal bounds \eqref{eq2.32}.  The detailed
arguments and related material appear in \cite{AL}. 

\subsection{The P\'{o}lya-Szeg\H{o} Conjecture for the Buckling of a
Clamped Plate} 

Much the same strategy as was used for the vibrating clamped plate 
can be used for the first eigenvalue $\Lambda_1(\Omega )$ of the 
buckling problem.  For $n=2$, this eigenvalue determines the critical 
buckling load of the clamped plate of the shape of $\Omega$ under 
uniform compressive loading around its boundary.  The analog of the 
Rayleigh conjecture in this setting,
\begin{equation}\label{eq2.43} \Lambda_1(\Omega )\geq
\Lambda_1(\Omega^\star)\quad \text{for } \Omega
\subset {\Bbb R}^n,\end{equation} 
with equality if and only if $\Omega$ is a ball, was conjectured 
by P\'{o}lya and Szeg\H{o} around 1950.  It is discussed in their 
book \cite{PS}.  At around the same time Szeg\H{o} proved this 
conjecture \cite{Sz1}, \cite{Sz3} under the assumption that 
$\Lambda_1(\Omega )$ has 
a first eigenfunction of fixed sign (which we can take to be 
nonnegative, i.e., $v_1\geq 0$  on $\Omega $).  Again this 
assumption turns out not to hold in all cases (see the article 
by Kozlov, Kondrat'ev, and Maz'ya \cite{KKM}).

The proof of the fixed-sign case (we assume $v_1 \geq 0$) can be 
accomplished much as was the corresponding case of the vibrating 
clamped plate above.  One simply starts from the variational principle
\begin{equation}\label{eq2.44} \Lambda_1(\Omega )=\min_{\varphi \in
H^2_0(\Omega )\atop \varphi
\not\equiv 0} \frac{\int_\Omega (\Delta\varphi )^2}{\int_\Omega |\nabla
\varphi |^2}=\frac{\int_\Omega (\Delta v_1)^2}{\int_\Omega |\nabla v_1|^2}
\end{equation} 
and proceeds as before.  One rearranges
$-\Delta v_1$ to obtain a symmetrized comparison problem.  Since
Talenti's version of the comparison argument (for (1.18) and 
(1.19)) can also be used to compare derivatives (indeed, this 
can be done directly from the integral-differential inequality for
$u^*$) we find easily that
\begin{equation}\label{eq2.45} \Lambda_1(\Omega )\geq
\frac{\int_{\Omega^\star}(\Delta \varphi )^2}{\int_{\Omega^\star}|\nabla
\varphi |^2}\end{equation} where $\varphi$ is now a radial function
on $\Omega^\star$ satisfying $\varphi (R)=0=\varphi'(R)$ (the proof
that $\varphi'(R)=0$ is exactly the same as that given for $v$ in
\eqref{eq2.35}).  Since this means that $\varphi$ is a valid trial
function for $\Lambda_1(\Omega^\star)$, the proof of
$\Lambda_1(\Omega )\geq \Lambda_1(\Omega^\star)$ is complete. 

For the general case (where $v_1$ is not necessarily of fixed sign)
we can also proceed as in the vibrating clamped plate problem.  That
is, we break $\Omega $ into the two parts $\Omega_\pm$ and rearrange 
$\mp \Delta v_1$ on each part separately to
obtain two symmetrized subproblems.  The hitch again comes at the
very last step, where the minimum over $a\in [0,1]$ turns out to
occur at $a^n=1/2$ for all dimensions $n\geq 2$.  Thus we do not
obtain  the general result \eqref{eq2.43} for any $n\geq 2$, but we
can again get nonoptimal lower bounds for
$\Lambda_1(\Omega)$ in the form
\begin{equation}\label{eq2.46} \Lambda_1(\Omega ) >
c_n\Lambda_1(\Omega^\star)\quad\text{for }
\Omega \subset {\Bbb R}^n,\end{equation} where the
dimension-dependent constants $c_n$ are found to tend to 1 as $n$
goes to infinity.  Some values of $c_n$ for low dimensions are
$c_2\approx 0.7877$, $c_3\approx 0.7759$, $c_4\approx 0.7872$,
$c_5\approx 0.8020$, $c_6\approx 0.8163$.  These results and related 
material can be found in \cite{AL} (see also \cite{ABL}).  

It is a rather remarkable fact that the bounds \eqref{eq2.46} with
precisely the same constants
$c_n$ can also be obtained by combining two classical eigenvalue
inequalities.  These are the inequality of Payne \cite{Pa1}
\begin{equation}\label{eq2.47} \Lambda_1(\Omega )\geq
\lambda_2(\Omega )\quad \text{for } \Omega
\subset {\Bbb R}^n,\end{equation} with equality if and only if
$\Omega $ is a ball, and Krahn's $\lambda_2$ inequality \cite{Kr2}
\begin{equation}\label{eq2.48} \lambda_2(\Omega )>
2^{2/n}\lambda_1(\Omega^\star)\quad \text{for }
\Omega \subset {\Bbb R}^n\end{equation} (this saturates as $\Omega $
disconnects into two equal disjoint balls).  Together these yield
\begin{equation}\begin{split}\label{eq2.49} \Lambda_1(\Omega
)>2^{2/n}\lambda_1(\Omega^\star) & = 2^{2/n} (C_n/|\Omega
|)^{2/n}j_{n/2-1,1}^2\\ &= 2^{2/n}\left(\frac{j_{n/2-1,1}}{j_{n/2
,1}}\right)^2 \left(\frac{C_n}{|\Omega |}\right)^{2/n}j_{n/2,1}^2\\
&= c_n\Lambda_1(\Omega^\star)\end{split}\end{equation} since
$\Lambda_1(\Omega^\star)=(C_n/|\Omega |)^{2/n}j_{n/2,1}^2$ and 
hence the constants $c_n$ are given by
\begin{equation}\label{eq2.50} c_n=2^{2/n}\left(
\frac{j_{n/2-1,1}}{j_{n/2,1}}\right)^2.\end{equation} The basic
observation here was made by Bramble and Payne \cite{BrPa}, though
they only gave the inequality for $n=2$ and hence did not
investigate the behavior of the $c_n$'s with varying $n$.  Neither
did they mention the connection between \eqref{eq2.49} and the
conjecture \eqref{eq2.43} of P\'{o}lya and Szeg\H{o}.  For further 
comments and details, one should consult \cite{AL} and/or \cite{ABL}.  

\section{Universal Inequalities for Eigenvalues}

\subsection{The General Inequalities of Payne-P\'{o}lya-Weinberger for
the Fixed Membrane Eigenvalues and their Extensions}

In this section we turn our attention to universal eigenvalue
inequalities.  We begin with a study of the Dirichlet eigenvalues of
$-\Delta $ on a bounded domain $\Omega \subset {\Bbb R}^n$.  This
subject began in 1955 with the work of Payne, P\'{o}lya, and
Weinberger \cite{PPW1}, \cite{PPW2}, who proved (among other things)
the bound
\begin{equation}\label{eq3.1} \lambda_{m+1}-\lambda _m\leq
\frac{2}{m}\sum^m_{i=1}
\lambda_i,\quad m=1,2,\ldots \end{equation} 
for $\Omega \subset {\Bbb R}^2$.  This result easily extends to 
$\Omega \subset {\Bbb R}^n$ as
\begin{equation}\label{eq3.2} \lambda_{m+1}-\lambda_{m}\leq
\frac{4}{mn}\sum^m_{i=1}\lambda_i,\quad m=1,2,\ldots .\end{equation}
Since Payne, P\'{o}lya, and Weinberger's paper \cite{PPW2},
\eqref{eq3.2} has been extended in several ways by a number of
authors.  For the Euclidean case there have been two main advances:
that of Hile and Protter \cite{HP} in 1980 and that of H.~C.\ Yang 
\cite{Yang} in 1991 (this paper is yet to be
published, as far as this author knows).  It turns out that, even
though the proofs given by Hile-Protter and Yang are dissimilar in 
several respects and are rather involved, the main approach is still
that of Payne-P\'{o}lya-Weinberger.  Moreover, we have succeeded in
streamlining these proofs to the point where, with the addition of
one key idea, they can be unified into a single overall approach
that is no harder than the original proof of \eqref{eq3.1} by Payne,
P\'{o}lya, and Weinberger.  The process of reducing these proofs to
their essence was begun in \cite{AB7},\cite{AB10} and is concluded in
\cite{Ash2}.  Since our approach is unified we head directly for the
result of H.~C.\ Yang, remarking only at the end on how the results of
PPW and Hile-Protter follow ``by simplification''.

To give an idea of Payne, P\'{o}lya, and Weinberger's basic method, we
begin by proving their bound
\begin{equation}\label{eq3.3} \frac{\lambda_2}{\lambda_1}\leq
1+\frac{4}{n} \quad \text{ for } \Omega \subset {\Bbb R}^n\end{equation} 
as a warm-up exercise (note that this is
the $m=1$ case of \eqref{eq3.2}).  The $n$-dimensional inequality 
\eqref{eq3.3} was first given explicitly by Thompson \cite{Thomp} 
in 1969, but certainly it (and more) is implicit in the work of 
Payne, P\'{o}lya, and Weinberger, as the generalization of their 
results to $\Omega \subset {\Bbb R}^n$ is entirely straightforward.  
To prove \eqref{eq3.3}, we introduce the trial function
\begin{equation}\label{eq3.4} \varphi =xu_1\end{equation} where $x$
represents any Cartesian coordinate $x_\ell (1\leq \ell \leq n)$ for
${\Bbb R}^n$.  By an appropriate choice of origin we can arrange that
$\int_\Omega x_\ell u^2_1=0$ for $1\leq \ell
\leq n$ (i.e., we put the origin at the center of mass of $\Omega $
for the mass distribution defined by $u^2_1$).  This guarantees
$\varphi \bot u_1$ (also $\varphi^{(\ell )}\bot u_1$ for
$1\leq \ell \leq n$ where $\varphi^{(\ell )}\equiv x_\ell u_1$) and
hence, by the Rayleigh-Ritz inequality (for $\lambda_2$),
\begin{equation}\label{eq3.5} \lambda_2\leq \frac{\int_\Omega
\varphi (-\Delta \varphi )}{\int_\Omega \varphi^2}.\end{equation}
Since
\begin{equation}\label{eq3.6} -\Delta \varphi =x(-\Delta
u_1)-2u_{1x}=\lambda_1xu_1-2u_{1x} =\lambda_1\varphi
-2u_{1x},\end{equation}
\eqref{eq3.5} gives
\begin{equation}\label{eq3.7} \lambda_2-\lambda_1\leq
-\frac{2\int_\Omega \varphi u_{1x}}{\int_\Omega
\varphi^2}.\end{equation} Now
\begin{equation}\label{eq3.8} 0\leq -2\int_\Omega \varphi
u_{1x}=-2\int_\Omega xu_1u_{1x}=-\int_\Omega x(u^2_1)_x=\int_\Omega
u^2_1=1\end{equation} where in the second-to-last step we integrated
by parts, and also by the Cauchy-Schwarz inequality,
\begin{equation}\label{eq3.9} (-2\int_\Omega \varphi u_{1x})^2\leq
4(\int_\Omega
\varphi^2)(\int_\Omega u^2_{1x}),\end{equation} implying (since
$-2\int_\Omega \varphi u_{1x}>0$)
\begin{equation}\label{eq3.10} \lambda_2-\lambda_1\leq
\frac{-2\int_\Omega \varphi u_{1x}}{\int_\Omega \varphi^2}\leq
\frac{4\int_\Omega u^2_{1x}}{-2\int_\Omega \varphi
u_{1x}}=4\int_\Omega u^2_{1x}.\end{equation} Obviously this same
argument applies to $\varphi^{(\ell )}\equiv x_\ell u_1$ for $1\leq
\ell \leq n$, allowing us to arrive at
\begin{equation}\label{eq3.11} \lambda_2-\lambda_1\leq 4\int_\Omega
u^2_{1x_\ell }\quad \text{for } 1\leq \ell \leq n.\end{equation} If
we now average these inequalities over $\ell $ we find 
\begin{equation}\label{eq3.12} \lambda_2-\lambda_1\leq
\frac{4}{n}\int_\Omega |\nabla
u_1|^2=\frac{4}{n}\lambda_1\end{equation} and hence \eqref{eq3.3}
follows.

Next we address the general results for $\lambda_{m+1}$.  The basic
strategy is the same as above in that we base our trial functions
$\varphi$ on $xu_i$ where $x$ is a Cartesian coordinate and
$u_i$ is a lower eigenfunction (i.e., $1\leq i\leq m$), but now to
enforce orthogonality to
$u_1,\ldots, u_m$ we no longer rely on the device of locating the
origin at the center of mass but rather we subtract away
counterterms which are just the projections of $xu_i$ along the
eigenfunctions $u_j$ for $1\leq j\leq m$ (when $m=1$, this amounts 
to shifting to the center of mass as before).  Thus as our trial
functions we take
\begin{equation}\label{eq3.13}
\varphi_i=xu_i-\sum^m_{j=1}a_{ij}u_j\quad \text{ for } 1\leq i\leq
m\end{equation} 
where
\begin{equation}\label{eq3.14} a_{ij}\equiv \int_\Omega
xu_iu_j=a_{ji}\end{equation} are the components of $xu_i$ along
$u_j$ for $1\leq j\leq m$ and thus clearly $\varphi_i\bot u_j$ for
$1\leq i,j\leq m$.  Also, it is straightforward to compute
\begin{equation}\label{eq3.15} \int_\Omega \varphi^2_i=\int_\Omega
\varphi _ixu_i=\int_\Omega
x^2u^2_i-\sum^m_{j=1}a^2_{ij}\end{equation} and
\begin{equation}\label{eq3.16}-\Delta
\varphi_i=\lambda_ixu_i-2u_{ix}-\sum
^m_{j=1}a_{ij}\lambda_ju_j\end{equation} so that (since $\varphi
_i\bot u_j$ for $1\leq j\leq m$)
\begin{equation}\label{eq3.17} \int_\Omega \varphi_i(-\Delta
\varphi_i)=\lambda_i\int_\Omega
\varphi^2_i-2\int_\Omega \varphi _iu_{ix}.\end{equation} It
therefore follows from the Rayleigh-Ritz inequality for
$\lambda_{m+1}$ that
\begin{equation}\label{eq3.18} (\lambda_{m+1}-\lambda_i)\int_\Omega
\varphi^2_i\leq -2\int_\Omega\varphi _iu_{ix}.\end{equation}
Proceeding much as before, we find
\begin{equation}\label{eq3.19}\begin{split} 0\leq -2\int_\Omega
\varphi_iu_{ix}&=-2\int_\Omega u_{ix}\left[
xu_i-\sum^m_{j=1}a_{ij}u_j\right]\\ &=-\int_\Omega
x(u^2_1)_x+2\sum^m_{j=1}a_{ij}\int_\Omega u_{ix}u_j\\
&=1+2\sum^m_{i=1}a_{ij}b_{ij}\end{split}\end{equation} 
where we have defined $b_{ij}$ by
\begin{equation}\label{eq3.20} b_{ij}=\int_\Omega
u_{ix}u_j.\end{equation} 
At this point, we can both recognize the $b_{ij}$'s as the 
components of the $u_{ix}$'s along the $u_j$'s and compute
them explicitly in terms of the $a_{ij}$'s.  It turns out that
these observations are both important for us.  (Also, $b_{ij}$ is
antisymmetric in $i$ and $j$, as can be recognized from
\eqref{eq3.20} immediately via integration by parts.)  We first
relate $b_{ij}$ to
$a_{ij}$:
\begin{equation}\begin{split}\label{eq3.21}
\lambda_ia_{ij}&=\int_\Omega (-\Delta u_i)xu_j\\ &=\int_\Omega
u_i[-\Delta (xu_j)]\\ &= \int_\Omega u_i[\lambda_jxu_j-2u_{jx}]\\
&=\lambda_ja_{ij}+2\int_\Omega u_{ix}u_j\\
&=\lambda_ja_{ij}+2b_{ij}\end{split}\end{equation} or
\begin{equation}\label{eq3.22}
2b_{ij}=(\lambda_i-\lambda_j)a_{ij}.\end{equation} Thus, from
\eqref{eq3.19},
\begin{equation}\label{eq3.23} 0\leq -2\int_\Omega
\varphi_iu_{ix}=1+\sum_{j=1}^m(\lambda_i-\lambda_j)a^2_{ij}.\end{equation}
Furthermore, by the Cauchy-Schwarz inequality,
\begin{equation}\label{eq3.24} \begin{split} (-2\int_\Omega \varphi_i
u_{ix})^2&=(-2\int_\Omega \varphi_i
[u_{ix}-\sum^m_{j=1}b_{ij}u_j])^2\\ &\leq 4(\int_\Omega
\varphi^2_i)(\int_\Omega [u_{ix}-\sum^m_{j=1}b_{ij}u_j]^2)\\
&=4(\int_\Omega \varphi^2_i)[\int_\Omega
u^2_{ix}-\sum^m_{j=1}b^2_{ij}].\end{split}\end{equation} From
\eqref{eq3.18} we now find
\begin{equation}\begin{split}\label{eq3.25}
(\lambda_{m+1}-\lambda_i)(\int_\Omega
\varphi^2_i)(-2\int_\Omega
\varphi_iu_{ix})&\leq (-2\int_\Omega \varphi_iu_{ix})^2\\ &\leq
4(\int_\Omega
\varphi_i^2)[\int_\Omega
u^2_{ix}-\sum^m_{j=1}b^2_{ij}]\end{split}\end{equation} and hence,
dividing by $\int_\Omega \varphi^2_i$,
\begin{equation}\label{eq3.26}
(\lambda_{m+1}-\lambda_i)(-2\int_\Omega \varphi_iu_{ix})\leq
4[\int_\Omega u^2_{ix}-\sum^m_{j=1}b^2_{ij}].\end{equation}
Inequality \eqref{eq3.26} holds even in the (unlikely) event that
$\varphi_i\equiv 0$, since in that case clearly its left-hand side
is 0 while its right-hand side is nonnegative (it is 4 times the
square of the norm of $u_{ix}-\sum^m_{j=1}b_{ij}u_j$).  Now, using
\eqref{eq3.22} and
\eqref{eq3.23}, \eqref{eq3.26} can be put in the form
\begin{equation}\label{eq3.27}
(\lambda_{m+1}-\lambda_i)[1+\sum^m_{j=1}(\lambda_i-\lambda_j)a^2_{ij}]
\leq 4\int_\Omega
u^2_{ix}-\sum^m_{j=1}(\lambda_i-\lambda_j)^2 a^2_{ij}\end{equation}
or
\begin{equation}\label{eq3.28}
\lambda_{m+1}-\lambda_i+\sum^m_{j=1}(\lambda_{m+1}-\lambda_j)
(\lambda_i-\lambda_j)a^2_{ij}\leq
4\int_\Omega u^2_{ix}.\end{equation}

From here it is clear how to finish our argument: we want to
introduce a factor of
$(\lambda_{m+1}-\lambda_i)$ to make the term involving $a^2_{ij}$
antisymmetric in $i$ and $j$, and then sum on $i$ from $1$ to $m$ to
make this term drop out.  We thus arrive at
\begin{equation}\label{eq3.29}\sum^m_{i=1}(\lambda_{m+1}-\lambda_i)^2\leq
4\sum^m_{i=1}(\lambda_{m+1}-\lambda_i)\int_\Omega
u^2_{ix}.\end{equation}  If we now recall that $x$ here could be any
$x_\ell (1\leq \ell \leq n)$ and sum over $\ell $ from $1$ to $n$ we
obtain
\begin{equation}\label{eq3.30}
n\sum^m_{i=1}(\lambda_{m+1}-\lambda_i)^2\leq 4\sum^m_{i=1}\lambda
_i(\lambda_{m+1}-\lambda_i)\end{equation} or
\begin{equation}\label{eq3.31} \sum^m_{i=1}(\lambda_{m+1}
-\lambda_i)\left(\lambda_{m+1}-\left( 1+\frac{4}{n}\right)
\lambda_i\right)\leq 0,\end{equation} which is the main (or
``first'') inequality of Hong Cang Yang \cite{Yang}.  This completes
our derivation of {\it Hong Cang Yang's first inequality}.

We now make a variety of comments about this inequality and the
inequalities of Payne-P\'{o}lya-Weinberger and Hile-Protter.  First,
noting that the left-hand side of
\eqref{eq3.31} is a quadratic in $\lambda_{m+1}$, we write it as
\begin{equation}\label{eq3.32} m\lambda^2_{m+1}-2\left(
1+\frac{2}{n}\right) \left(
\sum^m_{i=1}\lambda_i\right) \lambda_{m+1}+\left(
1+\frac{4}{n}\right)
\sum^m_{i=1}\lambda^2_i\leq 0\end{equation} and derive the explicit
upper bound
\begin{equation}\begin{split}\label{eq3.33} \lambda_{m+1}&\leq
[\text{ larger root of the quadratic }]\\  &=\frac{1}{m} \left[
\left( 1+\frac{2}{n}\right)\sum_{i=1}^m\lambda_i +\left\{ \left(
1+\frac{2}{n}
\right)^2 \left( \sum^m_{i=1}\lambda_i\right)^2-m \left(
1+\frac{4}{n}\right)
\sum^m_{i=1}\lambda^2_i\right\}^{1/2}\right].\end{split}\end{equation}
This bound is the best general upper bound yet derived by the
methods of Payne, P\'{o}lya, and Weinberger (or any other methods, for
that matter).  As was already observed by H.~C.\ Yang
\cite{Yang}, the bound \eqref{eq3.33} is much sharper than
previously known bounds for large $m$, since it comes much closer to
incorporating the Weyl asymptotic behavior of the eigenvalues
$\lambda_i$.  For further comments and observations in this regard,
see \cite{HaSt}, \cite{AB10},
\cite{Ash2}.

A simpler inequality due to H.~C.\ Yang, {\it Yang's second
inequality}, follows readily if we use the Cauchy-Schwarz inequality
to replace $m\sum^m_{i=1}\lambda^2_i$ by
$(\sum^m_{i=1}\lambda_i)^2$ on the right-hand side of \eqref{eq3.33}:
\begin{equation}\label{eq3.34} \lambda_{m+1}\leq \left(
1+\frac{4}{n}\right)
\frac{1}{m}\sum^m_{i=1}\lambda_i\quad \text{ for } m=1,2,\ldots
.\end{equation} This weaker inequality already implies the 
{\it Payne-P\'{o}lya-Weinberger inequality},
\eqref{eq3.2}, since we can obtain \eqref{eq3.34} from \eqref{eq3.2}
by replacing the $\lambda_m$ occurring on the left-hand side of
\eqref{eq3.2} by the average
$\frac{1}{m}\sum^m_{i=1}\lambda_i$, which is clearly less than
$\lambda_m$ (and is, in fact, strictly less for all $m\geq 2$ since
$\lambda_1<\lambda_i$ for all $i\geq 2$).

To obtain {\it Hile and Protter's inequality} we go back into our
proof of H.~C.\ Yang's inequality and make one modification.  In
inequality \eqref{eq3.24}, where we made use of the Cauchy-Schwarz
inequality, we use it, but without incorporating the counterterms
involving
$b_{ij}$'s, arriving at
\begin{equation}\label{eq3.35} (-2\int_\Omega \varphi_iu_{ix})^2\leq
4(\int_\Omega
\varphi^2_i)(\int_\Omega u^2_{ix}).\end{equation} 
This is how everyone had proceeded prior to H.~C.\ Yang.  No one had 
realized previously that one could make better (=``optimal'') use of 
the Cauchy-Schwarz inequality in this way (taking advantage of the 
known orthogonalities $\varphi_i\bot u_j$ for $1\leq i,j\leq m$ and the
fact that the $b_{ij}$'s have a simple relation to the
$a_{ij}$'s).  If our argument above is now carried through with
\eqref{eq3.35} replacing \eqref{eq3.24}, we arrive at
\begin{equation}\label{eq3.36}
(\lambda_{m+1}-\lambda_i)\left[1+\sum^m_{j=1}(\lambda_i-\lambda_j)
a^2_{ij}\right]
\leq 4\int_\Omega u^2_{ix}\end{equation} 
in place of \eqref{eq3.28}.

From here the clear thing to do to eliminate the unwanted (because
they're not easily controlled) terms in the $a^2_{ij}$'s is to 
divide through by $\lambda_{m+1}-\lambda_i$ and sum on $i$ from $1$ 
to $m$, yielding
\begin{equation}\label{eq3.37} \sum^m_{i=1}\frac{\int_\Omega
u^2_{ix}}{\lambda_{m+1}-\lambda_i}\geq \frac{m}{4}.\end{equation}
All the terms in $a^2_{ij}$ have dropped out due to antisymmetry.
Finally, recalling that $x$ could be any $x_\ell$, $1\leq \ell \leq
n$, and summing on $\ell$ we arrive at
\begin{equation}\label{eq3.38}
\sum^m_{i=1}\frac{\lambda_i}{\lambda_{m+1}-\lambda_i}\geq
\frac{mn}{4},\end{equation} 
which is the {\it Hile-Protter inequality}.  This inequality is stronger 
than the Payne-P\'{o}lya-Weinberger inequality \eqref{eq3.2}, since if we 
replace the $\lambda_i$ appearing in the denominator of the left-hand 
side of \eqref{eq3.38} by $\lambda_m$ we obtain
\eqref{eq3.2}.  It is also weaker than either of Yang's inequalities,
\eqref{eq3.33} or \eqref{eq3.34}.  It therefore follows that (for each 
$m=1,2,\ldots $)
\begin{equation}\label{eq3.39} \text{Yang 1}\Rightarrow \text{Yang 2}
\Rightarrow
\text{Hile-Protter}\Rightarrow \text{PPW}.\end{equation}

While our derivation above of the Hile-Protter inequality certainly
suggests that Yang's inequalities are stronger, it is not altogether
straightforward to show the middle implication in
\eqref{eq3.39}.  This is done in our longer paper \cite{Ash2}, which
is devoted to the topic of universal inequalities for the
eigenvalues of the Dirichlet Laplacian.  That paper also contains a
discussion of when the various inequalities given above are known to
hold strictly, as well as various other extensions and
generalizations.  A more complete set of references and some further
comments on them will be found there as well.

For anyone worried by our statement following inequality
\eqref{eq3.26} admitting the possibility that our trial functions
$\varphi_i$ might vanish identically and hence worried that our
derivation might result in a triviality, we make two remarks to
allay such fears.  The first is that our proof as given takes account
of these matters, and does indeed lead to a nontrivial inequality.
The second, and probably more useful, observation is that it can
never happen that all the $\varphi_i$'s vanish identically (if
indeed even one of them can vanish).  To see this directly, it
suffices to consider \eqref{eq3.23} and sum it on $i$ from 1 to $m$.
By antisymmetry the sum in $a^2_{ij}$ drops out, yielding
$m=-2\sum^m_{i=1}\int_\Omega \varphi_iu_{ix}$ and showing that not
all the $\varphi_{i}$'s can vanish identically (or even each be
orthogonal to
$u_{ix}$).  Thus at least some of the inequalities represented by
\eqref{eq3.26} will be nontrivial and will lead to a nontrivial
inequality \eqref{eq3.30} as the outcome of our final summations.

\subsection{Universal Inequalities for Low Fixed Membrane
Eigenvalues}

We next turn to some universal inequalities for low Dirichlet
eigenvalues, which again stem from the original work of Payne,
P\'{o}lya, and Weinberger \cite{PPW2}.  In 1956 they proved that for
$\Omega \subset {\Bbb R}^2$
\begin{equation}\label{eq3.40}
\frac{\lambda_2+\lambda_3}{\lambda_1}\leq 6.\end{equation} This
easily extends to $\Omega \subset {\Bbb R}^n$ as
\begin{equation}\label{eq3.41} \frac{\lambda_2+\lambda_3+\ldots
+\lambda_{n+1}}{\lambda_1}\leq
n\left(1+\frac{4}{n}\right)=n+4.\end{equation} These inequalities
are proved by using the Rayleigh-Ritz inequality in ``trace form''
or as usually applied, using rotations and translations to enforce
the further orthogonalities needed.  No $a_{ij}$'s appear (or can be
tolerated) in these arguments.  There are also a variety of
extensions of results of this type, the simplest being that of
Brands \cite{Bra} from 1964
\begin{equation}\label{eq3.42}
\frac{\lambda_2+\lambda_3}{\lambda_1}\leq
5+\frac{\lambda_1}{\lambda_2}\quad\text{ for } \Omega \subset {\Bbb
R}^2\end{equation} 
and its extension to ${\Bbb R}^n$
\begin{equation}\label{eq3.43} \frac{\lambda_2+\lambda_3+\ldots
+\lambda_{n+1}}{\lambda_1}\leq n+3+\frac{\lambda_1}{\lambda_2}\quad
\text{for } \Omega \subset {\Bbb R}^n\end{equation}
(see \cite{AB3} for a derivation, though this extension was certainly 
known to Hile and Protter earlier \cite{HP}).

Beyond this, much work has been done toward bounding the range of
values of $(\lambda_2/\lambda_1,\lambda_3/\lambda_1)$ for an arbitrary 
domain $\Omega \subset {\Bbb R}^2$ (with corresponding, but less 
extensive, work for $\Omega \subset {\Bbb R}^n$).  In particular, it 
would be very interesting to know the best bounds for 
$\lambda_3/\lambda_1$ and $(\lambda_2+\lambda_3)/\lambda_1$ (the 
best bound for $\lambda_2/\lambda_1$ is, of course, its value for 
a disk, given by \eqref{eq2.16}).  These quantities can be gotten at 
by looking at the range of possible values of
$(\lambda_2/\lambda_1,\lambda_3/\lambda_1)$, which is one 
motivation for its study.

The current state of knowledge regarding the range of values of
$(\lambda_2/\lambda_1,\lambda_3/\lambda_1)$ is summarized in the
paper of Ashbaugh and Benguria
\cite{AB10} (see also its precursors \cite{AB3} and \cite{AB6}).  In
particular, one should consult the graph given as Figure 1 on p.~38
of \cite{AB10}.  It is shown in \cite{AB10} that 
\begin{equation}\label{eq3.44} 5.077^+\leq \sup_\Omega
\frac{\lambda_2+\lambda_3}{\lambda_1}\leq 5.50661^-\end{equation}
and that
\begin{equation}\label{eq3.45} 3.1818^+\leq \sup_\Omega
\frac{\lambda_3}{\lambda_1}\leq 3.83103^-.\end{equation} The lower
bound in \eqref{eq3.44} is simply the value for a disk.  It is a
conjecture of Payne, P\'{o}lya, and Weinberger that this is the actual
maximum value of
$(\lambda_2+\lambda_3)/\lambda_1$ as well.  The lower bound in
\eqref{eq3.45} is $35/11$, the value taken by $\lambda_3/\lambda_1$
for a $\sqrt{8}$ by $\sqrt{3}$ rectangle (and is certainly the
maximum of $\lambda_3/\lambda_1$ among rectangles).  However, there
is no current guess as to the precise shape of domain that will
maximize $\lambda_3/\lambda_1$ (if a maximizer even exists).  The
best current thinking would have it be an elongated convex figure,
roughly in the shape of an oval or ellipse.  (Among ellipses the
largest value of $\lambda_3/\lambda_1$ seems to be very near to, but
slightly less than, 3.1818.)

\subsection{Universal Eigenvalue Inequalities in Other Spaces}

There are also versions of the Payne-P\'{o}lya-Weinberger inequality
\eqref{eq3.2} and its extensions for bounded domains in the constant
curvature spaces $S^n$ and $H^n$.  For example, for
$\Omega \subset S^2$ in 1975 Cheng \cite{Cheng} derived the bound
\begin{equation}\label{eq3.46} \frac{\lambda_2}{\lambda_1}\leq
1+2\left(\frac{2}{1+\cos
\Theta}\right)^4\quad \text{for } 0< \Theta <\pi\end{equation} where
$\Theta$ is defined as the outradius (=geodesic radius of the
circumscribing circle) of
$\Omega \subset S^2$.  This was improved to
\begin{equation}\label{eq3.47} \frac{\lambda_2}{\lambda_1}\leq
1+2\left( \frac{2}{1+\cos \Theta }\right) ^2\quad \text{for }
0<\Theta <\pi\end{equation} by Harrell \cite{Har} in 1993.  Both of
these results have extensions to higher eigenvalues and to
$S^n$.

On a slightly different front, by an extension of the method they
used to prove the Payne-P\'{o}lya-Weinberger conjecture for
$\lambda_2/\lambda_1$, Ashbaugh and Benguria proved
\cite{AB2}, \cite{AB7} (see Section 5) that for $\Omega\subset S^2$
\begin{equation}\label{eq3.48} \frac{\lambda_2}{\lambda_1}\leq
1+1.5387 \left(\frac{2}{1+\cos
\Theta }\right)^2\quad \text{for } 0<\Theta <\pi.\end{equation} 
The constant here is $j^2_{1,1}/j^2_{0,1}-1$ and comes from the proof 
of the Euclidean PPW conjecture as extended to general second-order 
elliptic operators in Section 4 of \cite{AB2}. 

Beyond this, there are the sharp PPW results for $\Omega $ contained 
in a hemisphere of $S^n$,
\begin{equation}\label{eq3.49} \lambda_2(\Omega )\leq
\lambda_2(B_{\lambda_1})\end{equation} where $B_{\lambda_1}$ is the
geodesic ball in $S^n$ having the same value of $\lambda_1$ as
$\Omega$, and
\begin{equation}\label{eq3.50} \frac{\lambda_2(\Omega
)}{\lambda_1(\Omega )}\leq
\frac{\lambda_2(\Omega^\star)}{\lambda_1(\Omega^\star)}\end{equation} 
where $\Omega^\star$ is the geodesic ball in $S^n$ having the same 
measure as $\Omega $.  These last two results are due to Ashbaugh and 
Benguria (see \cite{AB7}, \cite{AB11}, \cite{AB12}).  Note that 
\eqref{eq3.49} and \eqref{eq3.50} are much sharper than the bounds 
in terms of $\Theta$ listed previously, since the outradius of 
$\Omega $ is a much cruder measure of the size of $\Omega $ than 
either $\lambda_1(\Omega )$ or $|\Omega |$.  On the other hand, 
the bounds in terms of $\Theta$ apply to all domains $\Omega $, 
not just to those contained in a hemisphere.  While it may be 
possible that \eqref{eq3.49} and/or \eqref{eq3.50} hold beyond
the hemisphere, this is not proved as yet.  It might also be
mentioned that \eqref{eq3.50} follows easily from \eqref{eq3.49}
since it can be shown that $\lambda_2/\lambda_1$ for a geodesic ball
is an increasing function of its radius.  See \cite{AB7}, \cite{AB11},
and \cite{AB12} for more details and discussion.

Further results for domains $\Omega $ contained in $S^n$, $H^n$, or
other more general manifolds are due to Cheng \cite{Cheng} in 1975,
Li \cite{Li} in 1980, P.~C.\ Yang and S.-T.\ Yau
\cite{YangYau} in 1980, Lee \cite{Lee} and Leung \cite{Leung} in 1991,
Harrell \cite{Har} in 1993, and Harrell and Michel
\cite{HaMi1}, \cite{HaMi2} in 1994 and 1995.  Some of these papers
establish analogs of the Hile-Protter inequality in a more general
setting.  More discussion of this literature will be found in
\cite{Ash2}.  We mention, in particular, the Hile-Protter-type bound 
for a domain in $S^n$
\begin{equation}\label{eq3.51a}\sum^m_{i=1}\frac{4\lambda_i+n^2}
{\lambda_{m+1}-\lambda_i}\geq mn.\end{equation}
For $m=1$ this bound gives us 
\begin{equation}\label{3.51b} \lambda_2 \leq (1+\frac {4}{n}) \lambda_1+n.
\end{equation}
(Indeed, one has, more generally, the PPW-type bound 
\begin{equation}\label{3.51c} \lambda_{m+1} \leq \lambda_m + \frac{4}{mn}
\sum^m_{i=1} \lambda_i+n \leq(1+\frac{4}{n}) \lambda_m +n.)  \end{equation}
This bound can be regarded as an alternative and quite natural bound for
$\lambda_2$ in terms of $\lambda_1$ for domains $\Omega \subset S^n$ (to 
be compared with, say, \eqref{eq3.48}, at least when $n=2$).  Note that 
this bound ``takes account of'' the blow-up of $\lambda_2/\lambda_1$ 
that we expect when we approach the whole sphere $S^n$ by the presence 
of the constant term $n$ on its right-hand side.  Thus, even though 
$\lambda_1 \rightarrow 0$ as we approach the whole sphere while 
$\lambda_2 \rightarrow n$, this new bound gives us control of 
$\lambda_2$ without the need of a coefficient that blows up in this 
limit.  Moreover, the new bound $\lambda_2 \leq (1+\frac{4}{n})
\lambda_1 +n$ is seen to be sharp in this limit ({\it viz.}, for the full 
sphere $\lambda_2=n$, while the right-hand side becomes $n$ since 
$\lambda_1=0$).  The bound \eqref{eq3.51a} is actually the 
``HP-weakening'' of a stronger Yang-style bound for $\Omega 
\subset S^n$ first given in this context by H.~C.\ Yang \cite{Yang} 
(in the 1995 version of this paper).  For a full discussion of Yang's 
bounds and related inequalities (such as \eqref{eq3.51a}) in this 
context, see \cite{Ash2}.  Weaker, but generally more complex, precursors 
may be found in \cite{Cheng}, \cite{HaMi1}, \cite{Leung}, and 
\cite{YangYau} (most of these are directed at compact minimal 
hypersurfaces in $S^n$, but the ideas and calculations are much 
the same).      

\subsection{Universal Eigenvalue Inequalities for the Vibrating
Clamped Plate}

In their 1956 paper \cite{PPW2}, Payne, P\'{o}lya, and Weinberger also
established the bound
\begin{equation}\label{eq3.51} \Gamma_{m+1}-\Gamma_m\leq \frac{8}{m}
(\Gamma_1+\ldots +\Gamma_m),\quad m=1,2,\ldots ,\end{equation} 
for the eigenvalues $\Gamma_i$ of the clamped plate problem
(1.5) on a domain $\Omega \subset {\Bbb R}^2$.  This is the analog 
for the clamped plate of their bound \eqref{eq3.1} for the fixed 
membrane.  For $\Omega \subset {\Bbb R}^n$ this bound becomes
\begin{equation}\label{eq3.52} \Gamma_{m+1}-\Gamma_m\leq
\frac{8(n+2)}{n^2m}(\Gamma _1+\ldots +\Gamma _m),\quad m=1,2,\ldots
,\end{equation} 
and, even better,
\begin{equation}\label{eq3.53} \Gamma_{m+1}-\Gamma_{m}\leq
\frac{8(n+2)}{n^2m^2}(\Gamma_1^{1/2} +\ldots
+\Gamma_m^{1/2})^2,\quad m=1,2,\ldots .\end{equation} 
This latter inequality was not given by Payne, P\'{o}lya, and Weinberger, 
even for $n=2$, but it certainly could have been, as it is in a sense
implicit in their work.  Beyond this, in 1984 Hile and Yeh \cite{HY}, 
extending the approach of Hile and Protter to the clamped plate problem,
established the bound
\begin{equation}\label{eq3.54}
\sum^m_{i=1}\frac{\Gamma_i^{1/2}}{\Gamma_{m+1}-\Gamma_i}\geq
\frac{n^2m^{3/2}}{8(n+2)}\left( \sum^m_{i=1}\Gamma_i\right)
^{-1/2},\quad m=1,2,\ldots .\end{equation} Again, implicit in their
work is the better bound
\begin{equation}\label{eq3.55} \frac{n^2m^2}{8(n+2)}\leq
\left(\sum^m_{i=1}\frac{\Gamma^{1/2}_i}{\Gamma_{m+1}-\Gamma_i}\right)
\left(\sum^m_{i=1}\Gamma_i^{1/2}\right), \quad m=1,2,\ldots ,
\end{equation}
which was exhibited explicitly only later.  In fact Hook \cite{Hook}
in 1990 established \eqref{eq3.55} as a strict inequality.  Also in 1990, 
Chen and Qian \cite{ChQi} independently stated and proved \eqref{eq3.55}. 

Possibly \eqref{eq3.55} could be improved to
\begin{equation}\label{eq3.56} \left(
\sum^m_{i=1}\sqrt{\frac{\Gamma_i}{\Gamma_{m+1}-\Gamma_i}}\right)^2\geq
\frac{n^2m^2}{8(n+2)},\quad m=1,2,\ldots .\end{equation} 
This inequality would imply all the previous ones in this subsection
(\eqref{eq3.55}, for example, would follow easily using the Cauchy-Schwarz
inequality).  

It might be noted that the quantity $8(n+2)/n^2$ that appears in the 
foregoing inequalities is indeed natural from the PPW/HP point of view, 
just as $4/n$ is in the case of the fixed membrane.  This is because 
$8(n+2)/n^2$ arises as $(1+4/n)^2-1$ just as $4/n$ arises as $(1+4/n)-1$. 

Yet another inequality in this vein is   
\begin{equation}\label{eq3.57} \sum^m_{i=1}\frac{\Gamma_i}
{\Gamma_{m+1}-\Gamma_i}\geq\frac{n^2m}{8(n+2)},\quad m=1,2,\ldots .
\end{equation} 
This inequality is an easy consequence of \eqref{eq3.55} via Chebyshev's
inequality (it would also follow directly from \eqref{eq3.56} via the
Cauchy-Schwarz inequality).  Its chief appeal is its simplicity.  
It should be mentioned that as yet no one has
established an analog of H.~C.\ Yang's bound \eqref{eq3.31} (or
equivalently \eqref{eq3.33}) in this setting, i.e., for the
eigenvalues of the vibrating clamped plate.  So far \eqref{eq3.55}, 
which is best regarded as an analog of the Hile-Protter inequality 
\eqref{eq3.38} in this setting, is the closest anyone has come.

\subsection{Universal Eigenvalue Inequalities for the Buckling
Problem}

For the buckling problem for a clamped plate (problem (1.7))
even less is known, largely owing to the fact that the inner product
one must employ for this problem, $\langle f,g\rangle
\equiv \int_\Omega \nabla f\cdot \nabla g$, does not induce a
symmetric $a_{ij}$ matrix when one attempts the usual PPW approach
and sets $a_{ij}=\langle xv_i,v_j\rangle$.  This leads to all sorts
of extra complications, and thus far no one has been able to bring
the general case (for
$\Lambda_{m+1}-\Lambda_i$ or even just $\Lambda_{m+1}-\Lambda_m$) to
a satisfactory conclusion.  Payne, P\'{o}lya, and Weinberger
\cite{PPW2} were able to establish the low eigenvalue result
\begin{equation}\label{eq3.58} \Lambda_2/\Lambda_1<3\quad\text{for }
\Omega \subset {\Bbb R}^2.\end{equation} For $\Omega \subset {\Bbb
R}^n$ this reads
\begin{equation}\label{eq3.59}
\Lambda_2/\Lambda_1<1+4/n.\end{equation} Subsequently Hile and Yeh
\cite{HY} reconsidered this problem obtaining the improved bound
\begin{equation}\label{eq3.60} \Lambda_2/\Lambda_1\leq 2.5\quad
\text{for } \Omega \subset {\Bbb R}^2\end{equation} and, in general,
\begin{equation}\label{eq3.61} \frac{\Lambda_2}{\Lambda_1}\leq
\frac{n^2+8n+20}{(n+2)^2}\quad
\text{for } \Omega \subset {\Bbb R}^n.\end{equation} Both of these
inequalities can actually be shown to hold as strict inequalities.
It is of note that $\Lambda_2/\Lambda_1$ for a disk in 2 dimensions
is 1.796.  If the analog of the PPW conjecture held for this problem,
then this would be the best possible upper bound for
$\Lambda_2/\Lambda_1$ in 2 dimensions.

\subsection{More Inequalities for the Low Eigenvalues of the Clamped
Plate and Buckling Problems}

One can also derive inequalities analogous to \eqref{eq3.41} for the
clamped plate and buckling problems.  These read
\begin{align}\label{eq3.62} \frac{\Lambda_2+\ldots
+\Lambda_{n+1}}{\Lambda_1}&\leq n+4\quad
\text{ for } \Omega \subset {\Bbb R}^n,\\
\label{eq3.63} \frac{\Gamma ^{1/2}_2+\ldots + \Gamma
^{1/2}_{n+1}}{\Gamma^{1/2}_1}&\leq n+4\quad
\text{ for } \Omega \subset {\Bbb R}^n,\\
\intertext{and}
\label{eq3.64} \frac{\Gamma_2+\ldots +\Gamma_{n+1}}{\Gamma_1}&\leq
n+24\quad \text{ for } \Omega \subset {\Bbb R}^n.\end{align} 
Of the last two inequalities, \eqref{eq3.63} is the natural 
``PPW-analog'' and \eqref{eq3.64} is a weaker inequality that 
derives from it.  Note the presence of the PPW factor 
$\left( 1+\frac{4}{n}\right)$ in \eqref{eq3.62} and \eqref{eq3.63} 
(in the form $n+4=n\left(1+\frac{4}{n}\right)$).  The constant 
$n+24$ in \eqref{eq3.64} comes from the extreme case of \eqref{eq3.63} 
where $\Gamma_i^{1/2}/\Gamma_1^{1/2}=1$ for $2 \leq i \leq n$ 
and $\Gamma_{n+1}^{1/2}/\Gamma_1^{1/2}=5$.

Finally, we mention the known results for $\Gamma_2/\Gamma_1$.
Obviously (from \eqref{eq3.52}), Payne, P\'{o}lya, and Weinberger had 
\begin{equation}\label{eq3.65} \Gamma_2/\Gamma_1\leq (1+4/n)^2\quad
\text{ for } \Omega \subset {\Bbb R}^n\end{equation} 
(and even $\Gamma_{m+1}/\Gamma_m \leq (1+4/n)^2$ for all $m$).  
This was improved upon more recently by fairly elaborate means 
by Hile and Yeh \cite{HY}, who, for example, obtained the bounds 
$7.103$ and $4.792$ for dimensions $2$ and $3$ respectively. 
In general, these upper bounds are determined as the unique root larger
than $1$ of the cubic
\begin{equation}\label{eq3.66} (x-1)^3=\frac{512}{n^2(n+2)} x
\end{equation}
(though Hile and Yeh formulated their result in rather different terms).  
These values might be compared to those of the ball in $2$ and $3$ 
dimensions:  4.3311 and 3.2390, respectively.  Again the analog of the 
PPW conjecture for this problem would project these values as the 
optimal upper bounds for $\Gamma_2/\Gamma_1$ in $2$ and $3$ dimensions.
\\

As a final remark, we note that any result mentioned in this section 
(Section 3) without specific attribution is due to the author and that 
fuller details will be presented in forthcoming papers.

\section{Concluding Remarks and Open Problems}

In this concluding section we mention a few open problems and give some 
further hints to the literature. 

Of the original conjectures of Payne, P\'{o}lya, and Weinberger
\cite{PPW2}, only two remain open.  These are to show that
\begin{equation}\label{eq4.1}
\left.\frac{\lambda_2+\lambda_3}{\lambda_1}\leq
\frac{\lambda_2+\lambda_3}{\lambda_1}\right|_{\text{disk}}\approx
5.077 \text{ for } \Omega \subset {\Bbb R}^2,\end{equation} 
and the analogous result for $(\lambda_2+\ldots +\lambda_{n+1})/
\lambda_1$ for $\Omega \subset {\Bbb R}^n$, and that
\begin{equation}\label{eq4.2}
\left.\frac{\lambda_{m+1}}{\lambda_m}\leq
\frac{\lambda_2}{\lambda_1}\right|_{\text{ball}}
\quad \text{for all } \Omega \subset {\Bbb R}^n.\end{equation} 
The latter inequality is known for $m=1,2$, and $3$ (for $m=2$ and 
$3$ it follows from the stronger inequality $\lambda_4/\lambda_2\leq
\left.(\lambda_2/\lambda_1)\right|_{\text{ball}}$ proved in
\cite{AB5}), but as yet all higher cases remain open.  For further
discussion of these problems see \cite{AB3}, \cite{AB7}, and 
\cite{Ash1} (as well as \cite{AB5}).

There are also the well-known P\'{o}lya conjectures for the Dirichlet
and Neumann eigenvalues of the Laplacian on a domain $\Omega \subset 
{\Bbb R}^n$.  For the case of dimension 2, and with notation as in 
Section 1, these read
\begin{align}\label{eq4.3} \lambda_k & \geq \frac{4\pi
k}{A}\quad\text{for } k=1,2,\ldots\\
\intertext{and}
\label{eq4.4} \mu_k&\leq \frac{4\pi k}{A}\quad \text{for
}k=0,1,2,\ldots ,\end{align} 
where $A=|\Omega |$.  There are analogous conjectures for all 
dimensions $n>2$.  See \cite{AB7} for their statements and for 
further discussion.

For the vibrating clamped plate problem there remains the Rayleigh
conjecture
\begin{equation}\label{eq4.5} \Gamma_1(\Omega )\geq \Gamma_1
(\Omega^\star)\quad\text{for } \Omega
\subset {\Bbb R}^n\end{equation} for all $n\geq 4$, and, if one
could prove \eqref{eq4.5} (or in any event for $n=2,3$), the
PPW-type conjecture
\begin{equation}\label{eq4.6} \left.\frac{\Gamma_2}{\Gamma_1}\leq
\frac{\Gamma_2}{\Gamma_1}\right|_{\text{ball}} \quad\text{ for all
}\Omega \subset {\Bbb R}^n.\end{equation}

Similarly, for the buckling problem for the clamped plate there
remains the P\'{o}lya-Szeg\H{o} conjecture
\begin{equation}\label{eq4.7} \Lambda_1(\Omega )\geq
\Lambda_1(\Omega^\star)\quad\text{for } \Omega
\subset {\Bbb R}^n\end{equation} for all $n\geq 2$.  If this
conjecture could be proved then one could also consider the following
conjecture for ratios:
\begin{equation}\label{eq4.8} \left.\frac{\Lambda_2}{\Lambda_1}\leq
\frac{\Lambda_2}{\Lambda_1}\right|_{\text{ball}}\quad \text{for all
} \Omega \subset {\Bbb R}^n.\end{equation}

Obviously many other problems could be formulated and investigated.
For example, one could consider the ratios
$\Lambda_{m+1}/\Lambda_m$, $\Gamma_{m+1}/\Gamma_m$,
$(\Lambda_2+\Lambda_3)/\Lambda_1$, $(\Gamma_2+\Gamma_3)/\Gamma_1$,
and many other combinations analogous to those that have been considered 
for the eigenvalues of the Dirichlet Laplacian.  One could also consider
much of what has been discussed in this paper in the more general
setting of domains in Riemannian manifolds or for general second-order 
elliptic operators.

Other more extensive problem lists occur in the review papers of
Payne \cite{Pa2}, \cite{Pa3}, \cite{Pa4}, and in Yau's recent problem 
lists \cite{Yau1}, \cite{Yau2} (reprinted in \cite{SchYau}).  In addition, 
one could consult the final section of \cite{AB7} and also \cite{Ash1}.

\section*{Acknowledgements}

The author is grateful to Brian Davies and Yuri Safarov for the
opportunity to participate in the Instructional Conference on Spectral
Theory and Geometry in Edinburgh (March 29-April 9, 1998) and to the
International Centre for Mathematical Sciences for its generous
support.  In addition, he gratefully acknowledges his collaborators,
Rafael Benguria and Richard Laugesen, with whom a number of the
results summarized here were obtained.

\end{document}